\documentclass[10pt]{article}

\usepackage{latexsym,color,amsmath,amsthm,amssymb,amscd,amsfonts}

\newtheorem{theo}{Theorem}
\newtheorem{lem}[theo]{Lemma}
\newtheorem{cor}[theo]{Corollary}
\newtheorem{prop}[theo]{Proposition}
\newtheorem{defn}[theo]{Definition}

\def\R{\R}

\def\R{{\mathbb R}}
\def\grad{\nabla}

\def\eps{\varepsilon}

\def\qed{\hfill $\vcenter{\hrule height .3mm
\hbox {\vrule width .3mm height 2.1mm \kern 2mm \vrule width .3mm
height 2.1mm} \hrule height .3mm}$ \bigskip}

\def\eps{\varepsilon}

\def\lam{\lambda}

\def\to{\rightarrow}

\newcommand \supp{\operatorname{supp} \,}
\def\pmx{\begin{pmatrix}}
\def\emx{\end{pmatrix}}

\def\det{{\rm det}}
\def\supp{{\hbox{supp}}}
\def\Supp{{\hbox{supp}}}
\def\interior{{\hbox{int}}}

\def\R{\mathbb R}

\newcommand{\ps}{\psi^{\star}_{(s)}}

\begin{document}

\title{Divergence  for $s$-concave and log concave functions
\footnote{Keywords: entropy, divergence, affine isoperimetric inequalities, log Sobolev inequalities. 2010 Mathematics Subject Classification: 46B, 52A20,  60B}}

\author{Umut Caglar and Elisabeth M.  Werner\thanks{Partially supported by an  NSF grant}}

\date{}

\maketitle
\begin{abstract}
We prove new entropy inequalities for log concave and $s$-concave functions that
strengthen and generalize  recently established reverse log Sobolev  and Poincar\'e inequalities
for such functions.  This  leads naturally to the concept of $f$-divergence and,  in particular,  relative entropy   for 
$s$-concave and log concave functions. We establish their basic properties, among them the affine invariant valuation property.
Applications are given  in the theory of convex  bodies.
\end{abstract}

\section{Introduction}

There is a general approach to extend invariants and inequalities of convex bodies to the corresponding invariants and inequalities for functions.
Among the best known affine isoperimetric inequalities is the  Blaschke Santal\'o inequality \cite{Blaschke, MeyerPajor90, Santalo}.
The corresponding inequalities for log concave functions were proved by Ball \cite{KBallthesis} and Artstein, Klartag and Milman \cite{ArtKlarMil} (see also \cite{ MeyerFradelizi2007, Lehec2009}). 
A stronger inequality than the Blaschke Santal\'o inequality is the affine isoperimetric inequality  for convex bodies \cite{Blaschke,  Deicke, Santalo}.
The equivalent of this
 inequality   for log concave functions was established in \cite{ArtKlarSchuWer}: For  every 
log-concave
function $\varphi:\R^{n}\rightarrow [0, \infty)$ with enough smoothness and integrability properties  and such that $\int \varphi dx =1$, 
\begin{eqnarray}\label{equation1}
 \int_{\supp (\varphi)} \varphi \ \ln \bigg(\det \left( \text{Hess}\left(-\ln \varphi\right)\right)\bigg)  dx
 \leq 2 \left[  \operatorname{Ent}(\varphi) - \operatorname{Ent}(g) \right], 
\end{eqnarray}
where $g$ is the Gaussian,  
$\supp (\varphi)$ is the support,  
$\text{Hess}(\varphi) =  \left(\frac{\partial ^2 \varphi}{\partial x_i \partial x_j} \right)_{1 \leq i,j\leq n}$ is the Hessian 
 and
$ \operatorname{Ent}(\varphi)
=  \int_{\operatorname{supp}(\varphi)} \varphi \ln \varphi dx $
 is the entropy of $\varphi$. Thus,  the affine isoperimetric inequality  corresponds to a
reverse  log Sobolev inequality for
entropy. 
Equality holds in (\ref{equation1}) if and only if $\varphi(x)=C e^{-\langle A x, x \rangle}$, where
 $C$ is a positive constant  and $A$ is an $n \times n$  positive definite  matrix.  This characterization of  equality in  inequality  (\ref{equation1}) was achieved in \cite{CFGLSW}.
\par
Here, we strengthen and generalize inequality  (\ref{equation1}).  
\par
Inequality (\ref{equation1}) is  yet another instance of the 
rapidly developing,   fascinating connection 
between convex geometric analysis and information
theory. Further examples  can be found in  e.g., \cite{DCT, JenkinsonWerner, LutwakYangZhang2000, LutwakYangZhang2002/1, LutwakYangZhang2004/1, LutwakYangZhang2005,  PaourisWerner2011}.
In particular, it has been observed  \cite{Werner2012/1} that a fundamental notion of affine convex geometry, the $L_p$-affine surface area
can be viewed as  R\'enyi entropy from information theory, thus establishing a  link between  information theory and the 
powerful  $L_p$-Brunn-Minkowski
theory  \cite{Lutwak1996} of affine convex geometry.
Due to a number of highly influential works (see, e.g.,  
\cite{Ga3}-
\cite{HabSch2},  
\cite{Klain1996}, \cite{Klain1997},
\cite{Ludwig2002}-
\cite{LutwakZhang1997}, 
\cite{RuZ}, 
\cite{SW2002}- 
\cite{Werner2012/1},
\cite{WernerYe2008}, \cite{Zhang1994}), this theory 
is now  a central part of modern convex geometry. 
R\'enyi entropies are special cases of $f$-divergences whose  definition is given in Section \ref{Section-fdiv}. Such divergences and their related inequalities are  important tools  in  information
theory, statistics,  probability theory and machine learning \cite{BarronGyorfiMeulen, CoverThomas2006, GarciaWilliamson2012, HarremoesTopsoe, LieseVajda1987, LieseVajda2006, OsterrVajda, ReidWilliamson2011, WernerYe2013}. Consequently,  it is desirable to have such divergences available also in the theory of convex bodies   and this was achieved in  \cite{Werner2012}.
\par
In this paper, we further develop that direction.
We introduce $f$-divergences for  functions and establish some of their basic properties,  among them 
the affine invariance property and the valuation property. Valuations  were the critical ingredient in Dehn's solution of
Hilbert's third problem and, in the last decade,  have seen rapid growth as is demonstrated by  e.g.,   \cite{Alesker1999}-\cite{Alesker2003}, \cite{BernigFu}, \cite{Haberl2011}, \cite{Ludwig2002}-
\cite{Ludwig-Reitzner}, \cite{Schu}.
\par
We   prove the following  entropy inequality for log concave functions, i.e. functions 
of the form  $\varphi=e^{-\psi}$ with $\psi: \R^{n}\rightarrow \R$ convex. This inequality  is stronger than inequality (\ref{equation1}). Its proof 
 uses methods different from the ones used in \cite{ArtKlarSchuWer}.
\vskip 2mm
\begin{theo} \label{thm00}
 Let $f: (0, \infty) \rightarrow \mathbb{R}$ be a convex function.  Let $ \varphi:\R^{n}\rightarrow [0, \infty)$ be a
 log-concave
function.
Then
\begin{eqnarray}\label{thm00,1}
&&\hskip -14mm 
\int_{\Supp(\varphi)} \varphi  \  f \left( e^{\langle \frac{\nabla \varphi}{\varphi},x\rangle} \varphi^{-2} \left(  \det \left( Hess\left(-\ln \varphi\right)\right)\right)
\right) dx \nonumber \\ 
&&\hskip 20mm 
\geq \ f \left(  \frac{\int \varphi^\circ dx}{ \int \varphi dx}   \right)
 \  \left( \int_{\Supp(\varphi)} \varphi  dx \right).
 \end{eqnarray}
If  $f$ is concave, the inequality is reversed. 
If $f$ is linear, equality holds in (\ref{thm00,1}).  
Equality also holds 
if
 $\varphi(x)=C e^{-\langle A x, x \rangle}$, where
 $C$ is a positive constant  and $A$ is an $n \times n$ positive definite  matrix.
 \end{theo}
Here, $\nabla \varphi$ denotes the gradient of $\varphi$ and  $\varphi^\circ = \inf_{y \in \R^n } \left[ \frac{e^{- \langle x,y \rangle}}{\varphi(y)}  \right] $ \cite{ArtKlarMil} is the dual function of $\varphi$.
We will demonstrate  that the left hand side of the inequality (\ref{thm00,1}) is  {\bf the} natural definition of $f$-divergence $D_f(\varphi)$ for  a log concave
function $\varphi$,  so that inequality (\ref{thm00,1}) can be rewritten as 
\begin{equation}\label{thm00,2}
D_f(\varphi) \geq  \ f \left(  \frac{\int \varphi^\circ dx}{ \int \varphi dx}   \right)
 \  \left( \int_{\Supp(\varphi)} \varphi  dx \right).
 \end{equation}
This is shown in Section \ref{Section-logcon}. Inequality (\ref{thm00,2}) also holds for $s$-concave functions. We prove this  in Theorem \ref{s-thm}.
\par
If we let $f(t) = \ln t$ in Theorem  \ref{thm00},  we obtain the following corollary. 
\par
\begin{cor}\label{Korrolar01}
Let $ \varphi:\R^{n}\rightarrow [0, \infty)$ be a
log-concave
function. 
Then
\begin{eqnarray}\label{eqn0:Korollar1}  
\int_{\supp (\varphi)} \varphi \ \ln \bigg(\det \left( \text{Hess}\left(-\ln \varphi\right)\right)\bigg)  dx 
\leq 2  \operatorname{Ent}(\varphi) +  \|\varphi\|_{L^1} \ln\left[ e^n \left(\int \varphi \right)  \left(\int \varphi^\circ \right) \right], 
\end{eqnarray}
with equality if $\varphi(x)=C e^{-\langle A x, x \rangle}$, where
 $C$ is a positive constant  and $A$ is an $n \times n$  positive definite  matrix.
\end{cor}
\par
We show in Section \ref{Section-logcon} that 
inequality (\ref{eqn0:Korollar1}) involves the  {\em relative entropy} or {\em Kullback Leibler divergence} $D_{KL}(\varphi)$ (see 
Section \ref{Section-fdiv}  for the definition) of the function $\varphi$ and thus inequality (\ref{eqn0:Korollar1}) is equivalent to 
$$
D_{KL}(\varphi) 
 \leq \left( \int_{\Supp(\varphi)} \varphi  dx \right)   \ln\left(  \frac{\int \varphi^\circ dx}{ \int \varphi dx}   \right).
$$
Moreover, as it is  shown in Section \ref{Section-logcon},  the inequality of Corollary \ref{Korrolar01} is stronger than inequality (\ref{equation1}).
\par
It is important to note the affine invariant nature of the expressions (\ref{thm00,1}), (\ref{eqn0:Korollar1}) and of (\ref{eq:theo1}) below.
Both,  the respective left-hand sides and the right-hand sides,  are
invariant under volume-preserving linear transformations.
\par
The key ingredient to prove Theorem \ref{thm00} is (a special case) of  the following 
duality relation  for  log concave functions $\varphi:\R^{n}\rightarrow [0, \infty)$ and their  duals $\varphi^\circ$.
\par
\begin{theo}\label{theo1}
Let
$ \varphi:\R^{n}\rightarrow [0, \infty)$ be a log-concave
function. For a convex or concave function  $f: (0, \infty) \rightarrow \mathbb{R}$, let $f^*(t)=tf(1/t)$.
Then
\begin{equation} \label{eq:theo1}
D_f(\varphi^\circ) = D_{f^*}(\varphi).
\end{equation}
\end{theo}
We present several applications. 
In Section \ref{Section-specialf},    we consider  $f$-divergence for  special functions $f$, which,  on the level of convex bodies, correspond to  
$L_p$-affine surface areas. 
We refer to \cite{Lutwak1996, MW2, SW2004}
for the definition and to e.g.,   \cite{HabSch2}, \cite{Ludwig2010}, \cite{Ludwig-Reitzner1999}, \cite{Ludwig-Reitzner}, \cite{MW1}, \cite{SW1990}, \cite{SW2002}, \cite{Werner2007}-\cite{WernerYe2009} for more information on $L_p$-affine surface area for convex bodies. The $L_p$-affine surface areas for functions were  already introduced in \cite{CFGLSW}. Here, we establish several affine isoperimetric   inequalities  for these quantities. They are  the functional counterparts of known   inequalities for  convex bodies. 
Another application is given in Section \ref{Section-conbod}, where we apply our results about  log concave functions to convex bodies. Finally, in Section \ref{Section-lin}
we obtain  a reverse Poincar\'e inequality that is stronger than the one proved in \cite{ArtKlarSchuWer}.
\vskip 2mm
Throughout the paper we will assume that the convex or concave  functions $f: (0, \infty) \rightarrow \mathbb{R}$ and the 
$s$-concave and log concave  functions 
$ \varphi:\R^{n}\rightarrow [0, \infty)$ have enough smoothness and integrability properties so that the expressions 
considered in the statements make sense, i.e., we will always assume that $\varphi^\circ \in L^1(\supp (\varphi), dx)$, the Lebesgue integrable functions on the support of $\varphi$,  that
\begin{equation}\label{assume1}
\varphi \in C^2(\Supp(\varphi))  \cap L^1(\R^n, dx), 
\end{equation}
where  $C^2(\Supp(\varphi))$ denotes the twice continuously differentiable functions on their support, and that 
\begin{equation}\label{assume2}
 \varphi f 
\left(
\frac{e^{\frac{\langle \nabla \varphi, x\rangle}{\varphi}}
}
{\varphi^{2}} 
\mbox{det} \left(  \text{Hess} \left(-\ln \varphi \right) 
\right) \right) \in   L^1(\supp (\varphi), dx).
\end{equation}
See also  Remark (iv) after Definition \ref{defi2}.
\vskip 4mm

\section{ $f$-divergence for $s$-concave functions.} \label{Section-fdiv}
\subsection{Background on $f$-divergence.} \label{subsection:f-div}
In information theory, probability theory and statistics, an
$f$-divergence is a function that measures the difference between two (probability)
distributions. This notion was introduced by 
Csisz\'ar \cite{Csiszar}, and independently Morimoto \cite{Morimoto1963} and Ali \& Silvery \cite{AliSilvery1966}.
\par
Let $(X, \mu)$ be a measure space  and let  $P=p \mu$ and  $Q=q \mu$ be  (probability) measures on $X$ that are  absolutely continuous with respect to the measure $\mu$.  
Let $f: (0, \infty) \rightarrow  \mathbb{R}$ be a convex  or a concave  function.
The $*$-adjoint function $f^*:(0, \infty) \rightarrow  \mathbb{R}$ of $f$  is defined by 
\begin{equation}\label{adjoint}
f^*(t) = t f (1/t), \ \  t\in(0, \infty).
\end{equation}
\par
It is obvious   that $(f^*)^*=f$ and that $f^*$ is again convex  if $f$ is convex,  respectively concave if $f$ is concave.
Then the $f$-divergence   $D_f(P,Q)$ of the measures $P$ and $Q$ is defined by 
\begin{eqnarray}\label{def:fdiv1}
D_f(P,Q)&=&
 \int_{\{pq>0\} }f\left(\frac{p}{q} \right) q d\mu + f(0)\  Q\left(\{x\in X: p(x) =0\}\right)\nonumber \\
 &+& f^*(0) \ P\left(\{x\in X: q(x) =0\}\right),
\end{eqnarray}
provided the expressions exist. Here 
\begin{equation}\label{fat0}
f(0) = \lim_{t\downarrow 0} f(t)  \  \  \text{ and} \  \   f^*(0) = \lim_{t\downarrow 0} f^*(t).
\end{equation}
We make the convention that $0 \cdot \infty =0$. 
\vskip 2mm
Please note that 
\begin{equation}\label{fstern}
D_f(P,Q)=D_{f^*}(Q,P).
\end{equation}
With (\ref{fat0}) and as  
$$f^*(0) \ P\left(\{x\in X: q(x) =0\}\right) = \int _{\{q=0\}} f^*\left(\frac{q}{p} \right) p  d\mu =  \int _{\{q=0\}} f\left(\frac{p}{q} \right) q  d\mu,$$
we can write in short
\begin{equation}\label{def:fdiv2}
D_f(P,Q)=
 \int_{X} f\left(\frac{p}{q} \right) q d\mu.
\end{equation}
\vskip 2mm
Examples of $f$-divergences are as follows.
\par
\noindent
1.  For $f(t) = t \ln t$ (with  $*$-adjoint function $f^*(t) = - \ln t$), the $f$-divergence is   {\em Kullback-Leibler divergence} or {\em relative entropy} from $P$ to $Q$ (see \cite{CoverThomas2006})
\begin{equation}\label{relent}
 D_{KL}(P\|Q)= \int_{X} p \ln \frac{p}{q} d\mu.
\end{equation}
\par
\noindent
2. For   the convex or concave functions  $f(t) = t^\alpha$ we obtain the {\em Hellinger integrals} (e.g. \cite{LieseVajda2006})
\begin{equation}\label{Hellinger}
H_\alpha (P,Q) = \int _X  p^\alpha q^{1-\alpha} d\mu.
\end{equation}
Those are related to the 
R\'enyi divergence of order $\alpha$, $\alpha \neq 1$,  introduced by  R\'enyi \cite{Ren} (for $\alpha >0$) as 
\begin{equation}\label{renyi}
D_\alpha(P\|Q)=
\frac{1}{\alpha -1} \ln \left( \int_X p^\alpha q^{1-\alpha} d\mu \right)= \frac{1}{\alpha -1} \ln \left( H_\alpha (P,Q)\right).
\end{equation}
The case $\alpha =1$ is the relative entropy $ D_{KL}(P\|Q)$.

More on $f$-divergence can be found in e.g. \cite{GarciaWilliamson2012, LieseVajda1987,  LieseVajda2006, OsterrVajda, ReidWilliamson2011, Werner2012, WernerYe2013}.

\vskip 3mm
\subsection{ $f$-divergence for $s$-concave functions.}\label{subsection2.2}

Let $s \in \mathbb{R}$, $s \neq 0$. Let $\varphi: \mathbb{R}^n \rightarrow \mathbb{R}_+$.  Following Borell \cite{Borell1975}, 
we say that $\varphi$ is $s$-concave if for every $\lambda \in [0,1]$ and all  $x$ and $y$ such that $\varphi(x) >0$ and $ \varphi(y) > 0$,
\[
\varphi( (1-\lam)x + \lam y) \ge \left( (1-\lam) \varphi(x)^s + \lam \varphi(y)^s \right)^{1/s}.
\]
Note that $s$ can be negative. Now we want to define $f$-divergence for $s$-concave functions. To do that, let
$$P_\varphi ^{(s)}=  \frac{\mbox{det} \left[  \frac{-\text{Hess} \left( \varphi \right)}{\varphi}  + \left(1-s \right)
 \frac{\nabla \varphi \otimes \nabla \varphi
}{\varphi^2}  \right]} { \varphi \left(1- s \frac{\langle \nabla \varphi, x\rangle}{ \ \varphi}\right) ^{n+\frac{1}{s}} },  \hskip 3mm   Q_\varphi^{(s)}= \varphi  \left( 1 - s\frac{\langle \nabla \varphi, x\rangle}{ \varphi }\right).$$
Recall that we assume that the functions satisfy the conditions (\ref{assume1}) and (\ref{assume2}).
\vskip 2mm
\begin{defn}\label{defn:f-div-s-conc}
Let $f: (0, \infty) \rightarrow \mathbb{R}$ be  a  convex or  concave function.
Let $s \in \mathbb{R}$ and let $ \varphi:\R^{n}\rightarrow [0, \infty)$ be an $s$-concave function.
Then the $f$-divergence $D_f^{(s)}\left(P_\varphi^{(s)}, Q_\varphi^{(s)}\right)$ of $\varphi$ is 
\begin{eqnarray*}
D_f^{(s)}\left(P_\varphi^{(s)}, Q_\varphi^{(s)}\right)=
 \int_{\supp (\varphi)} \varphi\ f  \left(
\frac{\mbox{det} \left[  \frac{-\text{Hess} \left( \varphi \right)}{\varphi}  + \left(1-s \right)
 \frac{\nabla \varphi \otimes \nabla \varphi
}{\varphi^2}  \right]} { \varphi^2 \left(1- s \frac{\langle \nabla \varphi, x\rangle}{ \ \varphi}\right) ^{n+\frac{1}{s}+1} } \right)   \left( 1 - s\frac{\langle \nabla \varphi, x\rangle}{ \varphi }\right) dx.
\end{eqnarray*}
\end{defn}
\par
\noindent
We will sometimes  write in  short $D_f^{(s)}(\varphi)$ for $D_f^{(s)}\left(P_\varphi^{(s)}, Q_\varphi^{(s)}\right)$.
\par
\noindent 
Please note also that  for $s\neq1$  expression of the definition  can be rewritten as 
\begin{eqnarray}\label{sneq1}
&& D_f^{(s)}(\varphi)= D_f^{(s)}(P_\varphi^{(s)}, Q_\varphi^{(s)})=  \nonumber\\
 &&\int_{\supp (\varphi)} \varphi  f  \left(
\frac{\mbox{det} \left[  \left( \text{Hess} \left(-\ln \varphi \right) +
s \frac{\nabla \varphi \otimes \nabla \varphi
}{\varphi^2} \right) \right]} { \varphi^2 \left(1- s \frac{\langle \nabla \varphi, x\rangle}{  \ \varphi}\right) ^{n+\frac{1}{s}+1} } \right)   \left( 1 - s\frac{\langle \nabla \varphi, x\rangle}{ \varphi }\right) dx.
\end{eqnarray}
\par
\noindent
{\bf Remark.}
A similar expression holds for $D_{f}^{(s)}(Q_\varphi^{(s)}, P_\varphi^{(s)})$,  namely  
\begin{eqnarray}\label{stern}
&&\hskip -10mm D_{f}^{(s)}\left(Q_\varphi^{(s)}, P_\varphi^{(s)}\right) =  \nonumber \\
&&\hskip -13mm  \ \int_{\supp (\varphi)}  \  f  \left( 
\frac{ \varphi^2 \left(1- s \frac{\langle \nabla \varphi, x\rangle}{ \ \varphi}\right) ^{n+\frac{1}{s}+1} }{\mbox{det} \left[  \left(\text{Hess} \left(-\ln \varphi \right) +
s \frac{\nabla \varphi \otimes \nabla \varphi
}{\varphi^2} \right) \right]}  \right)  
 \frac{\mbox{det} \left[  \left(\text{Hess} \left(-\ln \varphi \right) +
s \frac{\nabla \varphi \otimes \nabla \varphi
}{\varphi^2} \right) \right]}{ \varphi \left( 1 - s \frac{ \langle \nabla \varphi, x\rangle}{  \varphi } \right)^{n+\frac{1}{s}}  }.
\end{eqnarray}
\par
\noindent
By (\ref{fstern}), $D_{f}^{(s)}\left(Q_\varphi^{(s)}, P_\varphi^{(s)}\right)=D_{f^*}^{(s)}\left(P_\varphi^{(s)}, Q_\varphi^{(s)}\right)= D_{f^*}^{(s)}(\varphi)$. Therefore it is enough to  only consider $D_f^{(s)}\left(P_\varphi^{(s)}, Q_\varphi^{(s)}\right)$. We will  do this throughout the paper.
\vskip 3mm
The motivation for this definition of $f$-divergence for $s$-concave functions comes from convex geometry.  In \cite{Werner2012}, $f$-divergence for a  convex body $K$ in $\mathbb{R}^n$ was introduced.
We refer to \cite{Werner2012} for more information and special cases and give here only the definition. 
\par
For $x \in \partial K$, the boundary of a sufficiently smooth  convex body $K$, let  $N_K(x)$ denote the outer unit normal to $\partial K$ in $x$ and let $\kappa_K(x)$ be the Gauss curvature in $x$.  $\mu_K$ is the usual surface area measure on $\partial K$. We put
\begin{equation}\label{densities}
p_K(x)=  \frac{\kappa_{K}(x)}{\langle x, N_K(x)\rangle^{n}} \, , \   \ q_K(x)=  \langle x, N_{K}(x) \rangle
\end{equation}
\noindent 
and 
\begin{equation}\label{PQ}
P_K=p_K\  \mu_K \ \ \ \text{and}   \ \ \   Q_K=q_K \ \mu_K.
\end{equation}
Then $P_K$ and $Q_K$  are   measures on $\partial K$ that are absolutely continuous with respect  to $\mu_{K}$. 
Let $f: (0, \infty) \rightarrow \mathbb{R}$ be a convex or concave function. In fact, $Q_K$ and $P_K$  are (up to the factor $n$) the {\em cone measures} (e.g.  \cite{PaourisWerner2011})  of $K$ and its polar 
\begin{equation}\label{polar-K}
K^\circ = \{y: \langle x, y \rangle \leq 1 \   \forall x \in K\},
\end{equation}
the latter  provided $K$ has sufficiently smooth boundary.
\par
 The 
$f$-divergence of $K$ with respect to the measures $P_K$ and $Q_K$ was defined in \cite{Werner2012} as 
\begin{eqnarray}\label{f-div1,0}
D_f(P_K, Q_K)=  \int _{\partial K} f\left(\frac{p_K}{q_K}\right) q_K d\mu_K 
= \int _{\partial K} f\left( \frac{  \kappa_K(x)} {\langle x, N_K(x) \rangle ^{n+1}}\right)  \langle x, N_K(x) \rangle 
d\mu_K.
\end{eqnarray}
\vskip 2mm
For $s >0$ such that $\frac{1}{s} \in \mathbb{N}$, we  associate with an $s$-concave function  $\varphi$ a  convex body $K_{s}(\varphi)$ \cite{ArtKlarMil} (see also \cite{ArtKlarSchuWer})
in $\R^{n}\times \R^\frac{1}{s}$, 
\begin{equation}\label{def.Ksf}
K_s(\varphi)=\big\{(x,y) \in \R^n \times \R^\frac{1}{s}: \sqrt {1/s} \  x \in
\overline{\supp (\varphi)}, \|y\| \leq
\varphi^s(\sqrt {1/s}\  x)\big\}.\end{equation}
The following proposition relates the definitions of $f$-divergence for the convex bodies and  $s$-concave functions. 
\vskip 2mm
\begin{prop} \label{motivation}
Let $s >0$ be such that $\frac{1}{s} \in \mathbb{N}$. Let $f: (0, \infty) \rightarrow \mathbb{R}$ be a convex or concave function. Let $ \varphi:\R^{n}\rightarrow [0, \infty)$ be an $s$-concave function. Then
\begin{equation}\label{div-funct1}
D_f^{(s)}\left(P_\varphi^{(s)}, Q_\varphi^{(s)}\right) = \frac{D_f\left( P_{K_s(\varphi)}, Q_{K_s(\varphi)}\right)}{s^\frac{n}{2} \mbox{vol}_{\frac{1}{s}-1}\left( S^{\frac{1}{s}-1} \right)}, 
\end{equation}
where $S^{\frac{1}{s}-1}$ is the $\left({\frac{1}{s}-1}\right)$-dimensional Euclidean sphere.
\end{prop}
\vskip 2mm
\begin{proof}
It was shown  in \cite{ArtKlarSchuWer} that for $z \in \partial (K_s(\varphi))$
\begin{eqnarray}\label{normale}
\langle z , N_{K_s(\varphi)} (z) \rangle = \frac{\varphi^s - \langle \grad\left(\varphi^s\right), x \rangle}{\left(1+ \|\grad \left(\varphi^s\right)\|^2 \right)^\frac{1}{2}}
\end{eqnarray}
and 
\begin{eqnarray}\label{expression}
&&\frac{ \kappa_{K_s(\varphi)}(z)} { \langle z , N_{K_s(\varphi)} (z) \rangle^{n+\frac{1}{s}+1} } = \frac{
\mbox{det} \left( \frac{ -\text{Hess} \  \varphi}{\varphi} +
\left(1-s \right) \frac{\nabla \varphi \otimes \nabla \varphi
}{\varphi^2} \right) }{\varphi^2\ \left(1- \sqrt{s}\frac{\langle \nabla \varphi, x \rangle}{ \varphi}\right)^{n+\frac{1}{s}+1}}, 
\end{eqnarray}
where $\varphi$ is evaluated at $\sqrt{1/s} x = (\sqrt{1/s}x_1,\ldots,\sqrt{1/s}x_n) \in \R^n$. 
We denote the collection of all
points $(x_1,\ldots,x_{n+\frac{1}{s}}) \in \partial K_s(\varphi)$ such that
$(\sqrt{1/s}x_1\ldots,\sqrt{1/s}x_n) \in \interior(\supp(\varphi))$ by $\tilde{\partial} K_s(\varphi)$. Since there is no contribution to the integral of $D_f\left(P_{K_s(\varphi)}, Q_{K_s(\varphi)}\right)$ from $\partial K_s(\varphi) \setminus \overline{\tilde{\partial} K_s(\varphi)}$
(since   the Gauss curvature vanishes on the part with full dimension, if it exists), 
we get 
with (\ref{normale}) and (\ref{expression}), 
\begin{eqnarray*}
&&D_f(P_{K_s(\varphi)}, Q_{K_s(\varphi)})= \int _{\partial K_s(\varphi)} f\left( \frac{  \kappa_{K_s(\varphi)}(z)} {\langle z, N_{K_s(\varphi)}(z) \rangle ^{n+1}}\right)  \langle z, N_{K_s(\varphi)}(z) \rangle 
d\mu_{K_s(\varphi)} \\
&&= \int_{\tilde{\partial} K_s(\varphi)} f\left(
\frac{
\mbox{det} \left( \frac{- \text{Hess} \  \varphi}{\varphi} +
\left(1-s \right) \frac{\nabla \varphi \otimes \nabla \varphi
}{\varphi^2} \right) }{\varphi^2\ \left(1- \sqrt{s} \frac{\langle \nabla \varphi, x \rangle}{ \varphi}\right)^{n+s+1}}\right)  
\frac{\left(\varphi^{s} -  \langle \nabla \left(\varphi^{s}\right) , x \rangle \right) } {  \left(1 + \| \nabla \left(\varphi^{s}\right) \|^2\right)^{\frac{1}{2}}  }  d\mu_{K_s(\varphi)} \nonumber \\
&&=  2\ \int_{\mathbb{R}^{n+\frac{1}{s}-1}}
f\left(
\frac{
\mbox{det} \left( \frac{- \text{Hess}   \varphi}{\varphi} +
\left(1-s \right) \frac{\nabla \varphi \otimes \nabla \varphi
}{\varphi^2} \right) }{\varphi^2\ \left(1- \sqrt{s} \frac{\langle \nabla \varphi, x \rangle}{  \varphi}\right)^{n+\frac{1}{s}+1}}\right)  
\frac{\varphi^{s} -  \langle \nabla \left(\varphi^{s}\right) , x \rangle  } {\varphi^{-s}  \left| x_{n+\frac{1}{s}}\right|} 
dx_1 \dots dx_{n+\frac{1}{s}-1}
\end{eqnarray*}
where $\varphi$ is evaluated at $\sqrt{1/s} x = (\sqrt{1/s}x_1,\ldots,\sqrt{1/s}x_n) $. The last
equality follows as the boundary of $K_s(\varphi)$ consists of two, ``positive'' and ``negative'',  parts.
As in \cite{ArtKlarSchuWer}, 
\begin{eqnarray*}
\int_{\R^{\frac{1}{s}-1}}\frac{ dx_{n+1} \dots dx_{n+\frac{1}{s}-1} }{
\left| x_{n+\frac{1}{s}}\right|} =
 \frac{\varphi^{1-2s}(\sqrt{1/s}x) }{2 s}  \   \mbox{vol}_{\frac{1}{s}}\left(B_2^{\frac{1}{s}}\right).
\end{eqnarray*}
Hence, 
\begin{eqnarray*}
&&D_f(P_{K_s(\varphi)}, Q_{K_s(\varphi)})= \\
&&c_s \int_{\left\{x: \sqrt{1/s}x  \in \supp (\varphi) \right\}} \varphi 
f\left(
\frac{
\mbox{det} \left( \frac{ -\text{Hess} \  \varphi}{\varphi} +
\left(1-s \right) \frac{\nabla \varphi \otimes \nabla \varphi
}{\varphi^2} \right) }{\varphi^2\ \left(1- \sqrt{s}\ \frac{\langle \nabla \varphi, x \rangle}{  \varphi}\right)^{n+\frac{1}{s}+1}}\right) 
\left(1- \sqrt{s} \frac{\langle \nabla \varphi, x \rangle}{  \varphi}\right) dx,
\end{eqnarray*}
where $\varphi$ is evaluated at $\sqrt{1/s} x = (\sqrt{1/s}x_1,\ldots,\sqrt{1/s}x_n) $ and  where $c_s = \frac{1}{s}\  \mbox{vol}_\frac{1}{s} \left(B^\frac{1}{s}_2\right)=
\mbox{vol}_{\frac{1}{s}-1}\left( S^{\frac{1}{s}-1} \right)$.
With the change of variable $\sqrt{1/s} x =y$, 
\begin{eqnarray*}
&&D_f(P_{K_s(\varphi)}, Q_{K_s(\varphi)})= \\
&& c_s s^{\frac{n}{2}} \int_{\supp (\varphi)} \varphi \ f  \left(
\frac{\mbox{det} \left[ - \left(\frac{\text{Hess} \left( \varphi \right)}{\varphi}  - \left(1-s \right)
 \frac{\nabla \varphi \otimes \nabla \varphi
}{\varphi^2} \right) \right]} { \varphi^2 \left(1- s \frac{\langle \nabla \varphi, y\rangle}{ \varphi}\right) ^{n+\frac{1}{s}+1} } \right) \  \left( 1 - s  \frac{\langle \nabla \varphi, y\rangle}{ \varphi }\right) dy= \\
&&c_s s^{\frac{n}{2}}  D_f^{(s)}\left(P_\varphi^{(s)}, Q_\varphi^{(s)}\right).
\end{eqnarray*} 
\end{proof}
\vskip 3mm
Now we describe some properties of the $f$-divergence for $s$-concave functions.
By (\ref{fstern}),  it is enough to  do  this for  $D_f^{(s)}(P_\varphi^{(s)}, Q_\varphi^{(s)})$ only.
\vskip 2mm
\begin{lem}\label{cor:f-div}
Let $ \varphi:\R^{n}\rightarrow [0, \infty)$ be an $s$-concave function and
let $f: (0, \infty) \rightarrow \mathbb{R}$ be a convex or concave function. Then 
  $D_f^{(s)}(\varphi)=D_f^{(s)}(P_\varphi^{(s)}, Q_\varphi^{(s)})$
is   invariant  under self adjoint $SL(n)$ invariant  linear maps and it is a valuation:
If $\max
(\varphi_1,\varphi_2)$ is $s$-concave,  then
\[ D_f^{(s)}(\varphi_1)+ D_f^{(s)}(\varphi_2) =  D_f^{(s)}(\max (\varphi_1, \varphi_2)) + D_f^{(s)}(\min (\varphi_1, \varphi_2)).\]
\end{lem}
\vskip 2mm
\begin{proof}[Proof]
Let $A: \mathbb{R}^n \rightarrow \mathbb{R}^n$ be a self adjoint,  $SL(n)$ invariant  linear map.
By Definition \ref{defn:f-div-s-conc}, 
\begin{eqnarray*}
&&\hskip -3mm D_f^{(s)}(P_{\varphi \circ A}^{(s)}, Q_{\varphi \circ A}^{(s)}) =
  \int_{\mbox{supp}( \varphi\circ A)}  \varphi(Ax)  \
\left (1- s \frac{ \langle \nabla \left(\varphi(Ax)\right), x \rangle}{ \  \varphi(Ax)}\right)\\
&& \hskip 35mm f \left(   \frac{
\mbox{det} \left(
\frac{ - \varphi(Ax) \  \mbox{Hess} \left( \varphi(Ax)\right) +  \left(1-s \right)\left(\nabla \left(\varphi(Ax)\right) \otimes \nabla \left(\varphi(Ax)\right) \right) }{   \left(\varphi(Ax)\right)^2}
\right)}{ (\varphi(Ax))^2 
\left (1- s \frac{ \langle \nabla \left(\varphi(Ax)\right), x \rangle}{ \  \varphi(Ax)}\right)^{n+\frac{1}{s}+1}} \right) \ dx \\ 
&&\hskip -3mm =  \frac{1}{|\det A|}  \int_{\supp (\varphi)}  \varphi   \left( 1 - s \frac{\langle \nabla \varphi, y\rangle}{ \varphi }\right)
 f  \left( (\det A)^2
\frac{\mbox{det} \left[ \frac{-\text{Hess} \left( \varphi \right)}{\varphi}  + \left(1-s \right)
 \frac{\nabla \varphi \otimes \nabla \varphi
}{\varphi^2}  \right]} { \varphi^2 \left(1- s \frac{\langle \nabla \varphi, y\rangle}{ \ \varphi}\right) ^{n+\frac{1}{s}+1} }    \right)    dy\\
&&\hskip -3mm= D_f^{(s)}(P_\varphi^{(s)}, Q_\varphi^{(s)}).
\end{eqnarray*}
\vskip 2mm
\noindent
\par
Next.  we establish 
the valuation property. There, $A^C$ denotes the complement of a set $A \subset \mathbb{R}^n$.
\begin{eqnarray*}
\hskip -7mm &&   D_f^{(s)}\left(P_{\varphi_1}^{(s)},  Q_{\varphi_1}^{(s)} \right)+ D_f^{(s)}(P_{\varphi_2}^{(s)},  Q_{\varphi_2}^{(s)} ) = \\
\hskip -7mm&& \int_{\supp (\varphi_1) \cap \supp(\varphi_2)} \varphi_1 f  \left(
\frac{\mbox{det} \left[  \frac{-\text{Hess} \left( \varphi_1 \right)}{\varphi_1}  + \left(1-s \right)
 \frac{\nabla \varphi_1 \otimes \nabla \varphi_1
}{\varphi_1^2}  \right]} { \varphi_1^2 \left(1- s \frac{\langle \nabla \varphi_1, x\rangle}{ \ \varphi_1}\right) ^{n+\frac{1}{s}+1} } \right)   \left( 1 - s\frac{\langle \nabla \varphi_1, x\rangle}{ \varphi_1 }\right) dx
+ \\
\hskip -7mm && \int_{\supp (\varphi_1) \cap (\supp (\varphi_2))^C } \varphi_1 f  \left(
\frac{\mbox{det} \left[  \frac{-\text{Hess} \left( \varphi_1 \right)}{\varphi_1}  + \left(1-s \right)
 \frac{\nabla \varphi_1 \otimes \nabla \varphi_1
}{\varphi_1^2}  \right]} { \varphi_1^2 \left(1- s \frac{\langle \nabla \varphi_1, x\rangle}{ \ \varphi_1}\right) ^{n+\frac{1}{s}+1} } \right)   \left( 1 - s\frac{\langle \nabla \varphi_1, x\rangle}{ \varphi_1 }\right) dx +\\
\hskip -7mm && \int_{\supp (\varphi_1) \cap \supp( \varphi_2)} \varphi_2 f  \left(
\frac{\mbox{det} \left[  \frac{-\text{Hess} \left( \varphi_2 \right)}{\varphi_2}  + \left(1-s \right)
 \frac{\nabla \varphi_2 \otimes \nabla \varphi_2
}{\varphi_2^2}  \right]} { \varphi_2^2 \left(1- s \frac{\langle \nabla \varphi_2, x\rangle}{ \ \varphi_2}\right) ^{n+\frac{1}{s}+1} } \right)  \left( 1 - s\frac{\langle \nabla \varphi_2, x\rangle}{ \varphi_2 }\right) dx
+ \\
\hskip -7mm && \int_{\supp (\varphi_2) \cap (\supp (\varphi_1))^C } \varphi_2 f  \left(
\frac{\mbox{det} \left[  \frac{-\text{Hess} \left( \varphi_2 \right)}{\varphi_2}  + \left(1-s \right)
 \frac{\nabla \varphi_2 \otimes \nabla \varphi_2
}{\varphi_2^2}  \right]} { \varphi_2^2 \left(1- s \frac{\langle \nabla \varphi_2, x\rangle}{ \ \varphi_2}\right) ^{n+\frac{1}{s}+1} } \right)  \left( 1 - s\frac{\langle \nabla \varphi_2, x\rangle}{ \varphi_2 }\right) dx\\
\hskip -7mm &&=  D_f^{(s)}(P_{\max (\varphi_1, \varphi_2)}^{(s)},  
Q_{\max (\varphi_1, \varphi_2)}^{(s)}) + D_f^{(s)}(P_{\min (\varphi_1, \varphi_2)}^{(s)},  
Q_{\min(\varphi_1, \varphi_2)}^{(s)}), 
\end{eqnarray*}
provided that $\max(\varphi_1, \varphi_2)$ is  $s$-concave.
\end{proof}
\vskip 3mm
Let $s \in \mathbb{R}$, $s \neq 0$  and let $ \varphi:  \mathbb{R}^n \rightarrow \mathbb{R}_+ $ be an $s$-concave function. 
 Let $\supp(\varphi)=\{x: \varphi(x) >0\}$ be the support of $\varphi$. Then $\supp(\varphi)$ is convex. We will assume throughout the rest of this section that $\supp(\varphi)$ is open and bounded, that $\varphi$ is $C^2$ on 
 $\supp(\varphi)$ and that  $\lim_{x \rightarrow \partial \supp(\varphi)} \varphi^s(x)=0$.
We   define a  function $\psi$ on $\supp(\varphi)$ (see also \cite{Rotem}) by
\begin{equation}\label{def: psi}
\psi(x)= \frac{1-\varphi^s(x)}{s},   \   \ x\in \supp(\varphi).
 \end{equation}
As $\varphi > 0$ on $\supp(\varphi)$, $\psi$ is well defined,  $\psi$ is convex on $\supp(\varphi)$, $ \psi < \frac{1}{s}$, if $s>0$ and $\psi > \frac{1}{s}$, if $s<0$.
We will use the following duality  definition from \cite{CFGLSW}.  First, 
 let $\left(\supp(\varphi)\right)^*= \{y: \sup_{x \in \supp(\varphi)} \langle x, y \rangle <1\}$.   
Note that  we can assume without loss of generality that $0 \in \supp(\varphi)$. 
If $0 \notin \supp(\varphi)$, pick $z \in \supp(\varphi)$ and consider 
 $\left(\supp(\varphi)-z \right)^* +z $.
Then $\left(\supp(\varphi)\right)^*$ is convex, open, bounded  and $0 \in \left(\supp(\varphi)\right)^*$. 
On the set $\frac{\left(\supp(\varphi)\right)^*}{s} $ we define 
\begin{equation}\label{def: psi-stern}
\ps (y) = \sup_{x \in \supp(\varphi)} \frac{\langle x, y \rangle - \psi(x)}{1 - s \psi(x)},   \  \  \\   y \in \frac{\left(\supp(\varphi)\right)^*}{s} .
\end{equation}
Then $\ps$ is convex,  and, as for $s >0$, $ \langle x, y \rangle < \frac{1}{s}$ for $ x \in  supp(\varphi)$ and $ y \in  \left(\supp(\varphi)\right)^*$,   we have that $ \ps< \frac{1}{s}$, if $s>0$ and, similarly,   that $\ps > \frac{1}{s}$, if $s<0$.
Observe also  that  for $s\rightarrow 0$ we obtain the Legendre transform $\mathcal{L} \psi(y) = \sup_{x} \left[\langle x, y \rangle - \psi(x)\right]$.
We denote 
$$\varphi_{(s)}^*(x) = \left(1- s \psi_{(s)}^*(x)\right)^{1/s}$$ 
the function corresponding to $\psi_{(s)}^*$.  $\varphi_{(s)}^*$ is well defined,  $s$-concave and, putting  $\varphi_{(s)}^* \equiv 0$ outside $\frac{\left(\supp(\varphi)\right)^*}{s} $,  coincides for $s >0$ with ${\cal L}_s(\varphi) (y) = \inf_{supp(\varphi)}  \frac{(1 - s \langle x, y \rangle)_+^{1/s}}{\varphi(x)}$ from \cite{ArtKlarMil}.  
\par
The  supremum  in (\ref{def: psi-stern}) is attained at $x$ such that
\begin{equation*}
\label{eq:attained}
y = \frac{1 - s \langle x, y \rangle}{1 - s \psi(x)} \, \nabla \psi(x)
\hbox{ which means } y = (1 - s \ps (y)) \nabla \psi (x). 
\end{equation*}
Moreover, 
\begin{equation}
\label{eq:gradient}
\frac{1}{1 - s \ps(y)} = \frac{1- s \psi(x)}{1 - s \langle x,y \rangle} = 1 + s ( \langle \nabla \psi(x), x \rangle - \psi(x)), 
\end{equation}
and  the relation between $y$ and $x$ is  
\begin{equation}
\label{eq:changevariable}
y = \frac{\nabla \psi(x)}{1 + s ( \langle \nabla \psi(x), x \rangle - \psi(x))}  = T_\psi (x).\
\end{equation}
It was noted in \cite{CFGLSW}   that the Jacobian is given by
\begin{equation*}
\label{eq:Jacobian}
dy = \left| \det \, dT_\psi(x) \right| dx =  \frac{1 - s \psi(x)}{\left(1 + s (\langle \nabla \psi(x), x \rangle - \psi(x))\right)^{n+1}} \ \det \, \text{Hess}\  \psi(x) \ dx.
\end{equation*}
It was also noted in  \cite{CFGLSW}  that the duality $(\ps)_{(s)}^{\star} = \psi$ holds and that therefore, 
\begin{equation}
\label{eq:JacobianDual}
\det \, \left(dT_{\psi}(x)\right) \det \,  \left( dT_{\ps}(y)\right)= 1.
\end{equation}
\par
Now, the next theorem provides a duality formula for an $s$-concave function $\varphi$ and $\varphi_{(s)}^*$.
It is  a generalization of a duality formula proved for special $f$ in \cite{CFGLSW}. 
\vskip 2mm
\begin{theo}\label{theo:duality-sconc}
Let  $f:(0,  \infty) \rightarrow \mathbb{R}$ be a convex or concave  function. Let   
$ \varphi:\R^{n}\rightarrow [0, \infty)$ be an $s$-concave function such that 
$C_\varphi$ is open and bounded, $\varphi$ is differentiable on 
 $C_\varphi$ and that  $\lim_{x \rightarrow \partial \supp(\varphi)} \varphi^s(x)=0$.
Then
\begin{equation}\label{s-polar-identity} 
 D^{(s)}_{f} (P^{(s)}_{\varphi_{(s)}^*},Q^{(s)}_{\varphi_{(s)}^*}) = D^{(s)}_{f^*}( P^{(s)}_{\varphi}, Q^{(s)}_{\varphi}).
\end{equation}
\end{theo}
\par
\noindent
{\bf Remark.} In particular, if $f\equiv 1$, or, equivalently, $f^*=Id$,  formula (\ref{s-polar-identity}) becomes 
\begin{eqnarray}\label{s-transformation}
(1+sn)   \int  \varphi_{(s)}^*
dx = \int_{\supp (\varphi)} 
\frac{\mbox{det} \left[  \frac{-\text{Hess} \left( \varphi \right)}{\varphi}  + \left(1-s \right)
 \frac{\nabla \varphi \otimes \nabla \varphi
}{\varphi^2}  \right]} { \varphi \left(1- s \frac{\langle \nabla \varphi, x\rangle}{ \ \varphi}\right) ^{n+\frac{1}{s}} }. 
\end{eqnarray}

\vskip 2mm
\begin{proof} 
By Definition \ref{defn:f-div-s-conc}, the change of variable (\ref{eq:changevariable}) and (\ref{eq:JacobianDual}) 
\begin{eqnarray*}
D^{(s)}_{f^*}\left( P^{(s)}_{\varphi}, Q^{(s)}_{\varphi}\right) &=& 
  \int_{\R^n} f^* \left( \frac{ \det \, \text{Hess}  \  \psi(x) }{ (1-s\psi(x))^{\frac{1}{s}-1}  | 1 - s \psi(x) + s  \langle x, \nabla \psi(x)\rangle  |^{\frac{1}{s} + n+1}} \right) \\
&& \hskip 1.8cm  \times \big(1 - s \psi(x) \big)^{\frac{1}{s} -1} | 1 - s (\psi(x) - \langle x, \nabla \psi(x)\rangle ) | \ dx \\
&=& \int_{\R^n} f^* \left( \frac{  \det \, dT_\psi(x) }{ (1-s\psi(x))^{\frac{1}{s}} \  | 1 - s \psi(x) + s  \langle x, \nabla \psi(x)\rangle  |^{\frac{1}{s}} } \right) \\
&& \hskip 3.3 cm  \times \big(1 - s \psi(x) \big)^{\frac{1}{s} -1}\big(1 - s \ps(y) \big)^{-1} \ dx \\
&=& \int_{\R^n} f \left( \frac{ (1-s\ps(y))^ {1-\frac{1}{s}}  \   \det \text{Hess} \   \ps(y) 
}{\left|1 - s \ps(y) + s  \langle y, \nabla \ps(y)\rangle  \right|^{\frac{1}{s} + n+1}
 }\right) \\
&& \hskip 0.6cm  \times \left(1 - s \ps(y) \right)^{\frac{1}{s}-1}   \left| 1 - s \ps(y) + s  \langle y, \nabla \ps(y)\rangle  \right| \ dy \\ 
&=& D^{(s)}_{f} \left(P^{(s)}_{\varphi_{(s)}^*},Q^{(s)}_{\varphi_{(s)}^*}\right).
\end{eqnarray*}
\end{proof}

The proof of the following   entropy inequality for $s$-concave functions is immediate with Jensen's inequality and identity (\ref{s-transformation}).
\vskip 2mm
\begin{theo} \label{s-thm}
 Let $f: (0, \infty) \rightarrow \mathbb{R}$ be a convex function.  Let $ \varphi:\R^{n}\rightarrow (0, \infty)$ be an $s$-concave function
 such that 
$C_\varphi$ is open and bounded, $\varphi$ is differentiable on 
 $C_\varphi$ and that  $\lim_{x \rightarrow \partial \supp(\varphi)} \varphi^s(x)=0$.
 Then
\begin{eqnarray*}\label{s-thm1}
D_f^{(s)}\left(P^{(s)}_{\varphi}, Q^{(s)}_\varphi\right) \geq \  (1+ns) \   \left( \int_{\Supp(\varphi)} \varphi  dx \right)   f \left(  \frac{\int  \varphi_{(s)}^*dx}{ \int \varphi dx}   \right).
 \end{eqnarray*}
If  $f$ is concave, the inequality is reversed. 
 \end{theo}

\vskip 4mm
\section{$f$-divergence  for log concave functions.}\label{Section-logcon}

A function $\varphi: \mathbb{R}^n \rightarrow [0, \infty)$ is log concave, if it is of the form $\varphi(x) = e^{-\psi(x)}$, where $\psi: \mathbb{R}^n \rightarrow \mathbb{R}$ is a convex function.
A  log-concave function $\varphi$ can be approximated
by the sequence of $k$-concave functions $ \{ \varphi_{k} \}_{k=1}^{\infty}$  
\begin{equation}\label{ApproxLogconcave}
\varphi_{k}=\left(1+ k \ln \varphi \right)_+^{\frac{1}{k}}, 
\hskip 10mm
k \in\mathbb N
\end{equation}
where  for $a \in \mathbb{R}$, $ a_+ = \max \{ a,0 \} $.  This motivates our definition for $f$-divergence for log concave functions.
We put 
\begin{equation}\label{Q,P}
Q_\varphi=\varphi  \hskip 4mm \text{and} \hskip 4mm   P_\varphi= \varphi^{-1}  e^{\frac{\langle\grad\varphi, x \rangle}{\varphi}} \mbox{det} \left[  \text{Hess} \left(-\ln \varphi \right)\right]
\end{equation}
and   define now the $f$-divergences for log concave  functions. 
\begin{defn} \label{defi2}
Let $f: (0, \infty) \rightarrow \mathbb{R}$ be  a  convex or  concave function and let $ \varphi:\R^{n}\rightarrow [0, \infty)$ be a log concave function.
Then the $f$-divergence $D_f(P\varphi, Q_\varphi)$ of $\varphi$ is 
\begin{equation}\label{div-Logconcave1}
D_f(P_\varphi, Q_\varphi)=
 \int_{\supp (\varphi)}\varphi  \ f  \left(
\frac{e^{\frac{\langle\grad\varphi, x \rangle}{\varphi}} }{\varphi ^2} \  \mbox{det} \left[  \text{Hess} \left( -\ln \varphi \right)\right] \right)dx.
\end{equation}
\end{defn}
\par
\noindent
Again, we will sometimes  write in  short $D_f(\varphi)$ for $D_f(P_\varphi, Q_\varphi)$.
\vskip 2mm
\noindent
{\bf Remarks and Examples.} 
(i) Similarly to (\ref{div-Logconcave1}), 
 \begin{eqnarray*}\label{div-Logconcave2}
 D_f(Q_\varphi, P_\varphi)=
  \int_{\supp (\varphi)}  \varphi^{-1} e^{\frac{\langle\grad\varphi, x \rangle}{\varphi}} \ \mbox{det} \left[ -  \text{Hess} \left(\ln \varphi \right)\right] \ f  \left( \frac{\varphi^2  }{
e^{\frac{\langle\grad\varphi, x \rangle}{\varphi}} \ \mbox{det} \left[ \text{Hess} \left(- \ln \varphi \right)\right]   }\right)dx.
\end{eqnarray*}
As by (\ref{fstern}), $D_f(Q_\varphi, P_\varphi)=D_{f^*}(P_\varphi, Q_\varphi)$, it is enough to consider $D_f(P_\varphi, Q_\varphi)$.
\par
\noindent
(ii) If we write a log concave function as $\varphi= e^{-\psi}$, $\psi$ convex, then (\ref{div-Logconcave1})
(and similarly  $D_f(Q_\varphi, P_\varphi)$) can be written as
\begin{equation}\label{div-Logconcave3}
D_f(P_\varphi, Q_\varphi)=
 \int_{\supp (\psi)} e^{-\psi}  \ f  \left(
e^{2 \psi - \langle\grad\psi, x \rangle}\  \mbox{det} \left[  \text{Hess} \psi \right] \right)dx.
\end{equation}
\par
\noindent 
(iii) Let $A$ be a positive definite, symmetric matrix, $C>0$ a constant and $\varphi(x) = C e^{- \langle Ax, x \rangle}$. Then
\begin{equation}\label{exponential2}
D_f(P_\varphi, Q_\varphi ) \ = \ f \bigg( \frac{ 2^n \det (A)}{C^2} \bigg) \frac{C \pi^{n/2}}{\sqrt{\det (A)}} .
\end{equation}
\par
\noindent 
(iv) Let $a$ be a non-zero vector in $\mathbb{R}^n$,  $C>0$ a constant  and  let $\varphi(x) = C e^{- \langle a, x \rangle}$. Then
\begin{equation*}\label{exponential1}
D_f(P_\varphi, Q_\varphi ) \ = \frac{C \ f (0)}{\prod_{i=1}^n a_i }  \left(\int _{\mathbb{R}}e^{-x} dx\right)^n, 
\end{equation*}
which is infinity, unless $f(0)=0$. Therefore, we require that $\varphi=e^{-\psi}$ is such that $\psi$ is strictly convex.

\vskip 2mm
If $\varphi$ is an $s_0$-concave function, then $\varphi$ is  $s$-concave for all $s \leq s_0$. In particular, $\varphi$ is log concave.
Thus $D_f(P_\varphi, Q_\varphi)$ is defined for $\varphi$ and   $D_f(P_\varphi, Q_\varphi) = D_f^{(0)}(P_\varphi^{(0)}, Q_\varphi^{(0)})$. 
\par
On the other hand, 
as it was remarked in  (\ref{ApproxLogconcave}),  every log  concave function can be  approximated by $s$-concave functions. The next Proposition shows that   Definition \ref{defi2} is compatible with  
Definition  \ref{defn:f-div-s-conc} for $s$-concave functions.
\vskip 2mm
\begin{prop} \label{prop-snach0}
Let $f : (0, \infty) \rightarrow \R$ be a convex or concave function.    Let $s >0$ and 
for a log concave function $ \varphi :\R^{n}\rightarrow [0, \infty) $ put $\varphi_s = \left(1+s \ln \varphi\right)_+^{\frac{1}{s}}$.
Then
\begin{equation*}
\lim_{s \rightarrow 0} D_f^{(s)}\left(P_{\varphi_s}^{(s)}, Q_{\varphi_s}^{(s)}\right)=D_f(P_\varphi, Q_\varphi).
\end{equation*}
\end{prop}
\vskip 2mm
\begin{proof}
Let $s>0$ and let $\varphi_s = \left(1+s \ln \varphi\right)_+^{\frac{1}{s}}$. Then
\begin{eqnarray*}
&&\hskip -7mm D_f^{(s)}\left(P_{\varphi_s}^{(s)},  Q_{\varphi_s}^{(s)}\right)= \int \left(1+s \ln \varphi\right)_+^{\frac{1}{s}}\left[1- \frac{s \langle \grad \varphi, x \rangle}{ \varphi 
\left(1+s \ln \varphi\right)_+} \right] \\
&& \hskip 35mm
f \left( 
\frac{\det
\left[ 
\frac{\mbox{Hess}(-\ln \varphi)}{\left(1+s \ln \varphi\right)_+^2} + \frac{s \nabla \left(\ln \varphi \right) \otimes  \nabla \left(\ln \varphi \right)}
{\left(1+s \ln \varphi\right)_+^2}  -  \frac{s \nabla \varphi  \otimes \nabla \varphi}{\varphi ^2\left(1+s \ln \varphi\right)_+^2}
 \right]}
{\left(1+s \ln \varphi\right)_+^\frac{2}{s} \left(1- \frac{s \langle \nabla \varphi, x \rangle}{ \varphi 
	\left(1+s \ln \varphi\right)_+ }\right) ^{n+\frac{1}{s} +1}}
\right) dx.
\end{eqnarray*}
Therefore
\begin{equation*}
\lim_{s \rightarrow 0} D_f^{(s)}\left(P_{\varphi_s}^{(s)}, Q_{\varphi_s}^{(s)}\right)=D_f(P_\varphi, Q_\varphi).
\end{equation*}
Note that we can interchange integration and limit because
conditions (\ref{assume1}) and (\ref{assume2}) hold. Compare also \cite{ArtKlarSchuWer}.
\end{proof}
\vskip 3mm
Similar to Lemma \ref{cor:f-div}, $f$-divergences for log concave functions are
affine invariant valuations.
Also, the proof is similar to the one of Lemma \ref{cor:f-div} and we omit it.
\vskip 2mm
\begin{cor}\label{cor:f-div-SLn}
Let $f: (0, \infty) \rightarrow \mathbb{R}$ be a convex or concave function and 
let $ \varphi:\R^{n}\rightarrow [0, \infty)$ be a 
log-concave function.  
Then $D_f(P_\varphi, Q_\varphi)$
is   invariant  under self adjoint $SL(n)$ maps and it is a valuation.
\end{cor}
\vskip 2mm
Recall that  for  $ \varphi :\R^{n}\rightarrow [0, \infty) $,  the dual function $\varphi^\circ$ \cite{ArtKlarMil}  is defined by  
$$ \ \varphi^\circ (x) = \inf_{y \in \R^n } \left[ \frac{e^{- \langle x,y \rangle}}{\varphi(y)}  \right].$$
This definition is connected with the Legendre transform $
\mathcal{L}\varphi(y) = \sup_{x \in \R^n} \left[ \langle x,y \rangle  -\varphi(x) \right] $, namely for $\varphi = e^{-\psi}$, 
\begin{equation}\label{Polare}
 \varphi^\circ = e^{-\mathcal{L}\left(- \ln \varphi \right) } = e^{-\mathcal{L}\left( \psi \right) }.
\end{equation}
\vskip 2mm
\noindent
{\bf Remark.}  Please observe that Proposition \ref{motivation}  justifies to 
call $Q_\varphi$ and $P_\varphi$ the {\em cone measures of the log-concave function $\varphi$ and its polar $\varphi^\circ$}.
\vskip 2mm
The next Theorem \ref{theo1}, already mentioned in the introduction, gives a duality relation  for a log concave function and its polar. 
We will see in Section \ref{Section-conbod} that it is the functional analogue of a duality formula for convex bodies.
\vskip 2mm
\noindent
{\bf Theorem 3.}
{\em Let  $f:(0,  \infty) \rightarrow \mathbb{R}$ be a convex or concave  function. Let
$ \varphi:\R^{n}\rightarrow [0, \infty)$ be a log-concave
function. Then
\begin{equation}\label{polar-identity} 
 D_{f}(P_{\varphi^\circ},Q_{\varphi^\circ}) = D_{f^*}( P_{\varphi}, Q_{\varphi}).
\end{equation}
}
\par
\noindent{\bf Remark.}  In particular, for $f\equiv 1$, or, equivalently, $f^*=Id$,  formula (\ref{polar-identity}) becomes 
\begin{eqnarray}\label{transformation}
 \int _{\Supp(\varphi)} \varphi^\circ
dx = \int_{\supp (\varphi)} \varphi^{-1}
\left(\mbox{det} \left(  \text{Hess} \left(-\ln \varphi \right) 
\right) \right) e^{\frac{\langle \nabla \varphi, x\rangle}{ \varphi}}. 
\end{eqnarray}
\begin{proof} 
We give a direct proof. But please observe that the proof also follows from Theorem \ref{theo:duality-sconc}, 
if we let $s \rightarrow 0$,
together with (\ref{Polare}) and Proposition \ref{prop-snach0}.
\par
We write  $ \varphi = e^{- \psi}$, $\psi$ convex,  and let  $\mathcal{L}  \psi (y)$ be the Legendre transform of $\psi$.
Please note that when $\psi$ is a $C^2$ strictly convex function, then 
$$
\psi(x) + \mathcal{L}  \psi(y) = \langle x, y \rangle
\hbox{ if and only if }
y = \nabla \psi(x)
\hbox{ if and only if }
x = \nabla \mathcal{L}  \psi(y).
$$
It follows that 
\begin{equation}
\label{dualitybis}
\forall y \in \R^n, \psi(\nabla \mathcal{L}  \psi(y)) = \langle y , \nabla \mathcal{L}  \psi (y) \rangle -  \mathcal{L}  \psi(y)
\end{equation}
and
\begin{equation}
\label{dualitygrad}
\nabla \psi \circ \nabla \mathcal{L}  \psi = \nabla \mathcal{L}  \psi \circ \nabla \psi = {\text Id}, 
\end{equation}
so that for any $x,y \in \R^n$,
\begin{equation}
\label{dualityhess}
\text{Hess} \, \psi (\nabla \mathcal{L}  \psi(y)) \ \text{Hess} \, \mathcal{L}  \psi(y) = {\text Id} = \text{Hess} \, \mathcal{L}  \psi(\nabla \psi(x)) \ \text{Hess} \,  \psi(x).
\end{equation}
Using equations (\ref{dualitybis}),  (\ref{dualitygrad}) and (\ref{dualityhess}), the change of variable $x=\nabla\mathcal{L} ( \psi (y))$ gives
\begin{eqnarray*}
 D_{f^*}( P_\varphi,Q_\varphi) &=& \int_{\mathbb R^n } \varphi  \  f^* \left(  \mbox{det} \left[\text{Hess} \left(- \ln  \varphi \right) \right] \frac{
 e^{\frac{\langle \nabla \varphi, x\rangle}{ \varphi}}  }{ \varphi ^2 }   \right) dx\\
&=& \int_{\R^n}   \det \left( \text{Hess} \psi (x)\right)  \  e^{ \psi(x) - \langle \nabla \psi, x\rangle}  f\left( \frac{e^{ -2 \psi(x) + \langle \nabla \psi, x\rangle} }{\det\left(\text{Hess} \psi (x)\right)  }  \right) dx \\
&=& \int_{\R^n} \det\left( \text{Hess} \psi ( \nabla\mathcal{L}  \psi (y))\right)   \  e^{ \psi ( \nabla\mathcal{L}  \psi (y)) - \langle y, \nabla (\mathcal{L}  \psi (y) ) \rangle }   \\
&& \hskip 1cm  \times  f \left( \frac{ e^{ -2 \psi( \nabla \mathcal{L}  \psi (y)) + \langle y, \nabla\mathcal{L}  \psi (y) \rangle} }{\det\left(\text{Hess} \psi ( \nabla \mathcal{L}  \psi (y) )\right)}\right) \   \det\left(\text{Hess}  \mathcal{L}  \psi (y) \right)  dy  \\
&=& \int_{\R^n}  e^{ - \mathcal{L}  \psi (y)  }   \  f \left( \det(\text{Hess} \mathcal{L}  \psi (y) ) \ e^{ - \langle y , \nabla \mathcal{L}  \psi  (y) \rangle +  2 \mathcal{L}  \psi (y) }  \right) \ dy  \\
&=& \int_{\mathbb R^n } \varphi^\circ  \ f \left(\mbox{det}\left[ \text{Hess} \left(-\ln  \varphi^\circ  \right) \right] \ \frac{
 e^{\frac{\langle \nabla \varphi^\circ, x\rangle}{ \varphi^\circ}}  }{ (\varphi^\circ) ^2 }   \right)  \\
&=& D_{f}(P_{\varphi^\circ},Q_{\varphi^\circ}).
\end{eqnarray*}
\end{proof}
\vskip 3mm
A consequence of Theorem \ref{theo1} is  the following  entropy inequality for log concave functions. This is Theorem \ref{thm00} of the introduction.
\vskip 2mm
\noindent
{\bf Theorem 1.}
{\em Let $f: (0, \infty) \rightarrow \mathbb{R}$ be a convex function and 
 let $ \varphi:\R^{n}\rightarrow [0, \infty)$ be 
 a log-concave
function.
Then
\begin{eqnarray*}\label{thm0,1}
D_f(P_\varphi, Q_\varphi)  \geq \ f \left(  \frac{\int \varphi^\circ dx}{ \int \varphi dx}   \right)
 \  \left( \int_{\Supp(\varphi)} \varphi  dx \right).
 \end{eqnarray*}
If  $f$ is concave, the inequality is reversed. 
If $f$ is linear, equality holds. Equality  also holds if 
 $\varphi(x)=C e^{-\langle A x, x \rangle}$, where
 $C$ is a positive constant  and $A$ is a $n \times n$  positive definite  matrix.}
\vskip 2mm
\begin{proof}
The inequality follows immediately from Jensen's inequality and identity (\ref{transformation}). 
Or, if we 
approximate $\varphi$ by $\varphi_s=(1+ s \ln \varphi)^\frac{1}{s}_+$, the inequality follows from Theorem \ref{s-thm} letting $s \rightarrow 0$.
\par
It is easy to check that equality holds if $f$ is linear and that  equality holds for $\varphi(x)=C e^{-\langle A x, x \rangle}$ by  (\ref{exponential2}). In fact, one can assume that $A$ is positive definite and symmetric.
\end{proof}
\vskip 3mm
If we let $f(t) = \ln t$ in Theorem  \ref{thm00},  we obtain the following corollary which is a reformulation of Corollary \ref{Korrolar01} of the introduction.
\vskip 2mm
\begin{cor}\label{Korrolar1}
Let $ \varphi:\R^{n}\rightarrow [0, \infty)$ be a 
log-concave
function.
Then
\begin{eqnarray}\label{eqn:Korollar1}
 D(P_\varphi || Q_\varphi) 
 \leq \left( \int_{\Supp(\varphi)} \varphi  dx \right)   \ln\left(  \frac{\int \varphi^\circ dx}{ \int \varphi dx}   \right).
\end{eqnarray}
Equality holds if  $\varphi(x)=C e^{-\langle A x, x \rangle}$, where
 $C$ is a positive constant  and $A$ is a $n \times n$  positive definite  matrix.
\end{cor}
\par
\noindent
{\bf Remarks.} 
\par
\noindent
(i) Inequality (\ref{eqn:Korollar1}) is stronger than (\ref{equation1}).
Indeed, as
\begin{eqnarray*} 
D(P_\varphi || Q_\varphi) = -2 \int  \varphi  \ln \varphi dx - n \int \varphi  dx + \int \varphi \ \ln \bigg(\det \left( \text{Hess}\left(-\ln \varphi\right)\right)\bigg)  dx, 
\end{eqnarray*}
inequality (\ref{eqn:Korollar1}) is equivalent to 
\begin{eqnarray*}
&&\hskip -5mm  \int_{\supp (\varphi)} \varphi \ \ln \bigg(\det \left( \text{Hess}\left(-\ln \varphi\right)\right)\bigg)  dx \nonumber \\
&&\hskip 5mm  \leq 2  \operatorname{Ent}(\varphi) +  \|\varphi\|_{L^1} \ln\left[ e^n 
\left(\int_{\supp (\varphi)} \varphi  dx\right)  \left(\int_{\supp (\varphi)} \varphi^\circ dx \right) \right].
\end{eqnarray*}
Now we apply the functional form of the Blaschke Santal\'o inequality \cite{ArtKlarMil,  KBallthesis, MeyerFradelizi2007, Lehec2009}. 
We assume without loss of generality  that $\int \varphi dx =1$.
Observe that we can also assume without loss of generality that $\int x \varphi(x) dx =0$. If  $\int x \varphi(x) dx =x_0$, replace $\varphi$ by $\tilde{\varphi}(x)=\varphi(x+x_0)$. We then get, 
\begin{eqnarray*}  \label{eqn2:Korollar1}
  \int_{\supp (\varphi)} \varphi \ \ln \bigg(\det \left( \text{Hess}\left(-\ln \varphi\right)\right)\bigg)  dx
 \leq 2  \operatorname{Ent}(\varphi) +  \ln\left( 2\pi e\right)^n \ ,
\end{eqnarray*}
which is inequality (\ref{equation1}).
\vskip 2mm
\noindent
(ii) The characterization of equality in (\ref{equation1}) now follows  by the equality characterization of the Blaschke Santal\'o inequality.
Indeed, from the arguments  in (i),  
if there is equality in (\ref{equation1}),  there is also
equality in the functional Blaschke-Santal\'o inequality. This implies that the function
has the form $\varphi(x)=C e^{-\langle A x, x \rangle}$  by  \cite{ArtKlarMil}.
\vskip 3mm
Let us state another corollary to  Theorem \ref{thm00}. Its proof follows immediately from Theorem \ref{thm00} and the functional 
Blaschke Santal\'o inequality and its equality characterization.
\vskip 2mm
\begin{cor} \label{cor-mono}
 Let $ \varphi:\R^{n}\rightarrow [0, \infty)$ be 
 a log-concave
function that  has center of mass at 0.
 Let $f: (0, \infty) \rightarrow \mathbb{R}$ be a convex, decreasing  function. Then 
\begin{eqnarray*}\label{thm0,2}
D_f\left(P_\varphi, Q_\varphi \right)  \geq   f \left(   \frac{(2 \pi)^n}{\left(\int \varphi dx \right)^2}   \right) \  \left( \int_{\Supp(\varphi)} \varphi  dx \right).
 \end{eqnarray*}
If  $f$ is a concave, increasing function, the inequality is reversed.
 \par
\noindent
Equality holds in both cases if and only if  $\varphi(x)= c  e^{-\langle A x, x \rangle}$, where $c$ is a positive constant and 
$A$ is an $n \times n$  positive definite  matrix.
\end{cor} 
\vskip 4mm

\section{Applications to special functions.} \label{Section-specialf}

Now we consider  special cases of $f$-divergences for log concave functions. Please recall that 
in subsection  \ref{subsection:f-div},  the   $\alpha$-R\'enyi entropies (\ref{renyi})  were introduced as  special $f$-divergences.
Examples  of such R\'enyi entropies  are,  for   log concave  functions $\varphi :\ R^{n}\rightarrow [0, \infty) $, for  $f(t)=t^\lambda$, $- \infty < \lambda <  \infty$, 
the {\em $L_\lambda$-affine surface areas $as_\lambda(\varphi)$} of $\varphi$, introduced in \cite{CFGLSW},  
\begin{equation}\label{asp-Logconcave}
as_\lambda(\varphi)=
\int_{\supp (\varphi)} \varphi \   \left(   
\frac{e^{\frac{\langle\grad\varphi, x \rangle}{\varphi}}}{\varphi^2}  \  \mbox{det} \left[   \text{Hess} \left(- \ln \varphi \right)\right] \right)^\lambda dx,
\end{equation}
or, writing $\varphi(x)= e^{-\psi(x)}$, $\psi$ convex,
\begin{equation}\label{sp-Logconcave1}
as_\lam( \varphi ) =as_\lam(e^{- \psi} ) = \int_{\R^n}e^{(2\lam-1)\psi(x)-\lam\langle x, \nabla\psi(x)\rangle}\left(\det \, {\text{Hess } \psi (x)}\right)^\lam dx.
\end{equation}
Especially,  $as_0(\varphi)= \int_{\supp (\varphi)} \varphi dx$ and, by (\ref{transformation}),  $as_{1}(\varphi) = \int_{\supp (\varphi)} \varphi^\circ dx$. 
Please note also  that for any log concave function  $\varphi $ we have that
$as_{\lambda} (\varphi) \geq 0 $.
Moreover,  by Corollary \ref{cor:f-div-SLn}, the  $as_{\lambda} (\varphi)$ are affine invariant valuations.
\vskip 3mm
We first want to give a definition for  $ as_{\infty} (\varphi)$ and  $ as_{-\infty} (\varphi)$,   similarly as it was done for convex bodies \cite{MW2}. To that end, for $\lambda >0$, 
let $ \tilde{as}_{\lambda } (\varphi) = \left( as_{\lambda} (\varphi) \right)^{ \frac{1}{\lambda}} \ $ and denote
$$ h =  \left(\varphi \right)^{ \frac{1}{\lambda}}   
\frac{e^{\frac{\langle\grad\varphi, x \rangle}{\varphi}}}{\varphi^2}    \mbox{det} \left[  \text{Hess} \left(- \ln \varphi \right)\right]  . $$
Denote by $\ \| f \|_{\lambda} = \left(\int  f^\lambda dx \right)^\frac{1}{\lambda}$ the $L_\lambda $ norm of a function $f$. Then, for $\lambda \rightarrow \infty$,

$$
\tilde{as}_\lambda(\varphi)= \left( 
\int_{\supp (\varphi)} h^\lambda dx. \right)^{\frac{1}{ \lambda}} \ = \  \| h \|_{\lambda}  \rightarrow \| h \|_{\infty} = \max_{x \in \supp (\varphi)}{ \frac{e^{\frac{\langle\grad\varphi, x \rangle}{\varphi}}}{\varphi^2}   \  \mbox{det} \left[   \text{Hess} \left(-\ln \varphi \right)\right]   }.
$$
Therefore, it is  natural to put
\begin{equation}\label{as+infinity}
as_{\infty} (\varphi) \ =   \  \max_{ x \in \supp (\varphi)}   \frac{e^{\frac{\langle\grad\varphi, x \rangle}{\varphi}}}{\varphi^2}  \  \mbox{det} \left[  \text{Hess} \left(-\ln \varphi \right)\right]  .        
\end{equation}
Similarly, for $\lambda \rightarrow - \infty$, 
\begin{eqnarray}\label{as-infinity}
\tilde{as}_\lambda(\varphi) &\rightarrow&   \frac{1}{ \max_{ x \in \supp (\varphi)}   \frac{e^{\frac{\langle\grad\varphi, x \rangle}{\varphi}}}{\varphi^2}  \  \mbox{det} \left[  \text{Hess} \left(-\ln \varphi \right)\right] } \nonumber \\
& = & \min_{ x \in \supp (\varphi)}  
\frac{\varphi^2}
{e^{\frac{\langle\grad\varphi, x \rangle}{\varphi}} \mbox{det} \left[  \text{Hess} \left(-\ln \varphi \right)\right] }  = as_{-\infty} (\varphi).        
\end{eqnarray}
Hence we have that 
$$
as_{-\infty} (\varphi) =\frac{1}{ as_{\infty} (\varphi)}.
$$
It is also easy to see that these expressions are  invariant under symmetric affine transformations with determinant $1$.
\vskip 3mm
The next theorem gives the analogue, for log concave functions, of a monotonicity behavior of  the $L_\lambda$-affine surface area that was established for 
convex bodies  in \cite{ Lutwak1996, WernerYe2008}. The case $\beta=0$ and $\alpha =1$ was already proved in \cite{CFGLSW}.
\vskip 2mm
\begin{prop} \label{theo:mono}
Let $\alpha \neq \beta,  \lambda \neq \beta$ be
real numbers. Let $\varphi: \R^{n}\rightarrow [0, \infty) $ be a log concave function.
\par
(i) If $1 \leq \frac{\alpha-\beta}{\lambda-\beta} < \infty$, then 
$
as_\lambda(\varphi)\leq \big(as_\alpha(\varphi) \big)^{\frac{\lambda-\beta}{\alpha-\beta}}
\big(as_\beta (\varphi)\big)^{\frac{\alpha-\lambda}{\alpha-\beta}}.
$
\par
(ii) If $1 \leq \frac{\alpha}{ \lambda} < \infty$, then 
$
as_\lambda(\varphi) \leq \left(as_\alpha(\varphi)\right)^\frac{\lambda}{\alpha} \left( \int \varphi \right)^\frac{\alpha-\lambda}{\alpha} .
$
\par
(iii)  If $ \beta \leq \lambda $, then 
$ as_\lambda(\varphi) \leq  \big(as_{ \infty }(\varphi) \big)^{ \lambda - \beta}
 \ as_\beta (\varphi).
$ 
 \vskip 2mm
\noindent
If $ \frac{\alpha-\beta}{\lambda-\beta}=1$ in (i),  respectively $\frac{\alpha}{\lambda}=1$ in (ii), then $\alpha=\lambda$ and equality holds trivially in (i)
respectively (ii). 
Equality also holds if $\varphi(x)= C e^{-\langle Ax,x\rangle}$.
\end{prop}
\vskip 2mm
\begin{proof}[Proof]
The proofs follow by  H\"older's inequality,  which, in (i), 
enforces the condition
$\frac{\alpha-\beta}{\lambda-\beta} > 1$.  The case(ii)  is a special case of (i) for $ \beta = 0$. We show  (i). The others follow similarly.
 \begin{eqnarray*} \label{firstcase}
 as_\lambda(\varphi)&=&\int_{\supp (\varphi)} \varphi \   \left(   
\frac{e^{\frac{\langle\grad\varphi, x \rangle}{\varphi}}}{\varphi^2}  \  \mbox{det} \left[  \text{Hess} \left( - \ln \varphi \right)\right] \right)^\lambda dx \\
&=&\int \left[\varphi \   \left(   
\frac{e^{\frac{\langle\grad\varphi, x \rangle}{\varphi}}}{\varphi^2}   \  \mbox{det} \left[  \text{Hess} \left( - \ln \varphi \right)\right] \right)^\alpha \ \right]
^\frac{\lambda-\beta}{\alpha-\beta} \\ 
 &&\hskip 15mm  \times \left[ \varphi \   \left(   
\frac{e^{\frac{\langle\grad\varphi, x \rangle}{\varphi}}}{\varphi^2}   \  \mbox{det} \left[  \text{Hess} \left( - \ln \varphi \right)\right] \right)^\beta \ \right]
^\frac{\alpha - \lambda}{\alpha - \beta} dx \\
&\leq & \big(as_\alpha(\varphi) \big)^{\frac{\lambda-\beta}{\alpha-\beta}}
\big(as_\beta (\varphi)\big)^{\frac{\alpha-\lambda}{\alpha-\beta}}.
 \end{eqnarray*}
\end{proof}
\vskip 3mm
\noindent
It follows from  Proposition \ref{theo:mono} (ii) that for $0 < \lambda \leq \alpha$,
$$
0 \leq \left(\frac{as_\lambda(\varphi) }{ \int \varphi dx }\right)^\frac{1}{\lambda} \leq \left( \frac{as_\alpha(\varphi) }{ \int \varphi dx } \right)^{\frac{ 1}{ \alpha}},
$$
which means that  for $\lambda >0$ the function 
$ \lambda \rightarrow\left( \frac{as_\lambda(\varphi) }{ \int \varphi}dx  \right)^\frac{1}{ \lambda} $  is bounded below by $0$ and
is  increasing for $ \lambda > 0$. 
Therefore, the  limit 
\begin{equation} \label{omega}
\Omega_{\varphi} =  \lim _{ \lambda \downarrow 0}  \left( \frac{as_\lambda(\varphi) }{ \int \varphi dx } \right)^\frac{1}{ \lambda} 
\end{equation} 
exists and the quantity $\Omega_\varphi$ is an  affine invariant. It is the analogue for log concave functions of an affine invariant   introduced by Paouris and Werner in \cite{PaourisWerner2011} for convex bodies. 
The quantity $\Omega_{\varphi}$ is related to the relative entropy as follows.
\vskip 2mm
\begin{prop}\label{affineinvariant-prop}
Let $\varphi: \R^{n}\rightarrow [0, \infty) $ be a log concave function. Then
$$\Omega_{\varphi} =  \lim _{ \lambda \downarrow 0}  \left( \frac{as_\lambda(\varphi) }{ \int \varphi dx } \right)^\frac{1}{ \lambda} 
=  \lim _{ \lambda \uparrow 0}  \left( \frac{as_\lambda(\varphi) }{ \int \varphi dx } \right)^\frac{1}{ \lambda}  = \exp 
\left(  \frac{ D( P_{\varphi}  | |  Q_{\varphi}) }{\int \varphi dx } \right).$$
\end{prop}
\vskip 2mm
\begin{proof}[Proof]
By definition and de l'H\^ospital, 
\begin{eqnarray*}
 \Omega_{\varphi}&=&  \lim _{ \lambda \downarrow 0}  \left( \frac{as_\lambda(\varphi) }{ \int \varphi} \right)^\frac{1}{ \lambda}   
 = \lim _{ \lambda \downarrow 0} \  \exp \left( \frac{1}{\lambda} 
   \ln \left( \frac{ as_\lambda(\varphi) }{ \int \varphi dx } \right) 
 \right)  \\ 
 &=& \exp \left(  \lim _{ \lambda \downarrow 0}  \frac{\int \frac{d}{d\lambda} \bigg[ \varphi  \  \left(\frac{e^{\frac{\langle\grad\varphi, x \rangle}{\varphi}}}{\varphi^2}  \  \mbox{det} \left[  \text{Hess} \left(-\ln \varphi \right)\right]  \right)^\lambda \bigg] dx}{as_\lambda(\varphi)} 
 \right) \\
&=& \exp \left(   \frac{\int \varphi \  \ln \left(\frac{e^{\frac{\langle\grad\varphi, x \rangle}{\varphi}}}{\varphi^2}  \  \mbox{det} \left[  \text{Hess} \left(-\ln \varphi \right)\right] dx \right)}{ \int \varphi  dx } 
 \right) \\
&=& \exp \left(  \frac{ D(P_\varphi | | Q_\varphi) }{ \int \varphi dx}\right).
\end{eqnarray*}
It also follows from Proposition \ref{theo:mono} (ii) that for $ \lambda <0$, the function $ \lambda \rightarrow\left( \frac{as_\lambda(\varphi) }{ \int \varphi}dx  \right)^\frac{1}{ \lambda} $  is  increasing 
We  compute $ \lim _{ \lambda \uparrow 0}  \left( \frac{as_\lambda(\varphi) }{ \int \varphi^\circ} \right)^\frac{1}{ \lambda} $ as above.
\end{proof}
\vskip 3mm
\begin{cor} \label{affineinvariant-cor}
Let $ \varphi :\R^{n}\rightarrow [0, \infty) $ be a log concave function. 
\par
(i) $\Omega_{\varphi}  \leq    \left( \frac{as_\lambda(\varphi) }{ \int \varphi dx } \right)^\frac{1}{ \lambda}$ for all $  \lambda > 0$ and  
$\Omega_{\varphi}  \geq     \left( \frac{as_\lambda(\varphi) }{ \int \varphi dx } \right)^\frac{1}{ \lambda}$ for all $  \lambda < 0 $.
\par
(ii)
$ \Omega_{\varphi}  \ \Omega_{\varphi^\circ}  \ \leq  \ 1. $
\par
(iii) 
$ \Omega_{\varphi} = \lim_{\alpha \rightarrow 1}  \left( \frac{as_\alpha(\varphi^\circ ) }{ \int \varphi dx } \right)^\frac{1}{ 1- \alpha}$.
\par
\noindent
Equality holds in (i) and (ii)  if $ \varphi =  C e^{-\langle Ax,x\rangle} $.
\end{cor}
\vskip 2mm
\begin{proof}[Proof]
(i) is deduced immediately from the monotonicity behavior of the function $ \lambda \rightarrow\left( \frac{as_\lambda(\varphi) }{ \int \varphi}dx  \right)^\frac{1}{ \lambda} $ and the definition of  $\Omega_\varphi$.
\par
\noindent
(ii) By (i),  $ \Omega_{\varphi}  \leq   \frac{as_1 (\varphi) }{\int \varphi dx }= \frac{\int\varphi  ^\circ dx }{\int\varphi dx  }$ and  
$ \Omega_{\varphi^\circ}  \leq   \frac{as_1 (\varphi^\circ) }{\int \varphi ^\circ dx }= \frac{\int\varphi  dx }{\int\varphi ^\circ dx  }$.
Here, we have also used  the bipolar property $(\varphi^\circ)^\circ = \varphi$. Thus  (ii) follows.
\par
\noindent
(iii)  We use the  duality formula $as_\lambda (\varphi) = as_{ 1 -\lambda} (\varphi^\circ) $ which was first proved in \cite{CFGLSW}. Note that it can also be obtained as a special case of Theorem \ref{theo1} for $f(t)=t^\lambda$.
By definition
\begin{eqnarray*}
 \Omega_{\varphi^\circ} &=&  \lim _{ \lambda \rightarrow 0}  \left( \frac{as_\lambda(\varphi)^\circ }{ \int \varphi^\circ} \right)^\frac{1}{ \lambda} 
 = \lim _{ \lambda \rightarrow 0}  \left( \frac{as_{1-\lambda}(\varphi) }{ \int \varphi^\circ} \right)^\frac{1}{ \lambda}
 = \lim _{ \alpha \rightarrow 1}  \left( \frac{as_{\alpha}(\varphi) }{ \int \varphi^\circ} \right)^\frac{1}{1- \alpha}.
\end{eqnarray*}
Therefore, $ \Omega_{\varphi} = \lim_{\alpha \rightarrow 1}  \left( \frac{as_\alpha(\varphi^\circ ) }{ \int \varphi dx } \right)^\frac{1}{ 1- \alpha}$.
\end{proof}

\vskip 4mm

\section{Application to convex bodies.}\label{applications} \label{Section-conbod}
Let us now consider the case of 2-homogeneous functions $\psi$, that is $\psi(\lambda x) = \lambda^2 \psi(x)$ for any $\lambda \in \R_+$ and $x \in \R^n$. Such functions $\psi$ are necessarily (and this is obviously sufficient) of the form $\psi(x) = \|x\|_K^2/2$ for a certain convex body $K$ with $0$ in its interior, where we have denoted by $\| . \|_K$ the gauge function of $K$,
\begin{eqnarray*}
\|x\|_K = \min\{\lambda \geq 0: \lambda x \in K\} = \max_{y \in K^\circ}\langle x, y \rangle = h_{K^\circ}(x).
\end{eqnarray*}
Differentiating with respect to $\lam$ at $\lam=1$, we get  
\begin{equation}\label{grad=psi}
\langle x, \nabla \psi(x) \rangle = 2 \psi(x). 
\end{equation} 
Now we apply this function to the identities and inequalities which we have obtained for $f$-divergences for log concave functions. It was already observed in \cite{CFGLSW}
that the $L_\lambda$-affine surface area for log concave functions is a generalization of $L_\lambda$-affine surface area for convex bodies, Indeed, it was noted there that  if one applies  the log concave function $ \varphi_K= \exp \left( -\frac{\|\cdot\|_K^2}{2} \right) $ to  Definition \ref{asp-Logconcave}, then one obtains $L_\lambda$-affine surface area for convex bodies.
Please note also that 
\begin{equation}\label{integralphi}
 \int e^{ -\frac{\| x \|_K^2}{2} } \ dx = \frac{(2 \pi)^{ \frac{n}{2}} |K|}{|B^n_2|}   \ \ \text{and} \ \  \int e^{ -\frac{\| x \|_{K^\circ}^2}{2} } \ dx = \frac{(2 \pi)^{ \frac{n}{2}} |K^\circ |}{|B^n_2|} .
\end{equation}
Recall that $B^n_2$ denotes the $n$-dimensional Euclidean unit ball, and for a convex body $K$ in $\mathbb{R}^n$, $K^\circ$ is the polar of $K$ (\ref{polar-K}) 
and $|K|$ is its volume.
\vskip 3mm
The following is a generalization of a   result in \cite{CFGLSW} but it  is proved  in exactly the same  way. We include the proof for completeness.
\vskip 2mm
\begin{theo} Let $K$ be a convex body in $\mathbb{R}^n$ with $0$ in its interior. Let $\varphi_K = \exp \left( -\frac{\|\cdot\|_K^2}{2} \right)$ and              $f:( 0, \infty) \rightarrow  \mathbb{R}$ be a convex or concave function. Then
\begin{eqnarray*}
D_f(P_{\varphi_K}, Q_{\varphi_K}) 
= \frac{(2\pi)^\frac{n}{2}}{n|B_2^n|} \ D_f(P_K, Q_K ).
\end{eqnarray*}
Here, $P_K$ and $Q_K$ are as in (\ref{PQ}) and, for $\varphi_K$, $P_{\varphi_K}$ and $Q_{\varphi_K}$ are as in (\ref{Q,P}).
\end{theo}
\vskip 2mm
\begin{proof}
We will use  formula (\ref{div-Logconcave3}) for $\psi=\frac{\|\cdot\|_K^2}{2}$ and integrate in polar coordinates with respect to the normalized cone measure  $\frac{P_K}{n|K|}$  (\ref{PQ}) of  $K$. Thus, if we write $x=r z$, with $z \in\partial K$, $dx=r^{n-1}dr dP_K(z)$. We also use that  the map  $x\mapsto\det \, {\text{Hess } \psi (x)}$ is $0$-homogeneous. Therefore,  with (\ref{grad=psi}), 
\begin{eqnarray*}
 D_f(P_{\varphi_{K}}, Q_{\varphi_K}) &=& \int_0^{+\infty}r^{n-1}e^\frac{-r^2}{2}dr \int_{\partial K} f \left(\det \, {\text{Hess} \, \psi (z)}\right) \ dP_K(z)\\
 &=& \frac{2\pi)^\frac{n}{2}}{n |B_2^n|}\int_{\partial K}  f \left(\det \, {\text{Hess} \, \psi (z)}\right) \ dP_K(z).
\end{eqnarray*}
It was proved in  \cite{CFGLSW}  that  for all $z\in\partial K$, 
\begin{equation}\label{hesse}
\det  \,( \mathrm{Hess}_z \psi ) = \frac{ \kappa_K (z) }{\langle z,  N_K(z)\rangle ^{n+1} }.
\end{equation}
 Observe  that for the $ G_K(z)$ introduced in this lemma,  $\Vert G_K(z) \Vert_{K^\circ} =\langle z,N_K(z)\rangle$. Thus 
\begin{eqnarray*}
D_f(P_{\varphi_K}, Q_{\varphi_K})&=& \frac{(2\pi)^\frac{n}{2}}{n|B_2^n|}\int_{\partial K}  f \left(\frac{ \kappa (x) }{ \langle x,N_K(x)\rangle^{n+1} }\right) \langle x,N_K(x)\rangle d\mu_K(x)\\
&=& \frac{(2\pi)^\frac{n}{2}}{n|B_2^n|} D_f(P_K, Q_K ).
\end{eqnarray*}
\end{proof}
\vskip 2mm
Now we apply Theorem \ref{thm00} to $\varphi = \exp \left( -\frac{\|\cdot\|_K^2}{2} \right)$ and  we obtain the following inequalities. Those were already proved, with different methods,  in \cite{Werner2012}. In fact, it was shown there that equality holds if and only if $K$ is an ellipsoid.
\vskip 2mm
\begin{cor}\label{Korrolar2}
Let $K$ be  a convex body in $\mathbb{R}^n$ with the origin in its interior. 
Let $f:( 0, \infty) \rightarrow  \mathbb{R}$ be a convex function. Then 
\begin{equation}\label{pre-asa1}
D_f(P_K, Q_K) \geq 
 n \ |K| f \left( \frac{|K^\circ|}{|K|}\right).
\end{equation}
If $f$ is concave, the inequality is reversed.
Equality holds 
if $K$ is an ellipsoid.
\end{cor}
\vskip 2mm
\begin{proof} Let $\varphi = \exp \left( -\frac{\|\cdot\|_K^2}{2} \right) $ and let $f$ be convex.   By
Theorem \ref{thm00}, together with (\ref{integralphi}),  
\begin{eqnarray*}
&& \int_{\Supp(\varphi)} e^{ -\frac{\| x \|_K^2}{2} } \  f \left(  \left[  \det \left( Hess\left(\frac{\| x \|_K^2}{2} \right)\right)\right] 
\right)dx
  \geq \ 
  \frac{(2 \pi)^{ \frac{n}{2}} |K|}{|B^n_2|}  f\left(  \frac{|K^\circ|}{ |K|}   \right).
 \end{eqnarray*}
Now we  use again (\ref{hesse}) 
and, as above,  make the change of variable,  $x=r z$, $z \in \partial K$. Then 
this becomes 
\begin{eqnarray*}
\int_{\partial K}  f \left( \frac{ \kappa (x) }{ \langle x,N_K(x)\rangle^{n+1} }\right)\langle x,N_K(x)\rangle d\mu_K(x) 
  \geq n\   |K| f\left( \frac{|K^\circ|}{|K|}\right).
\end{eqnarray*}
For $f$  concave,  the direction in the inequality changes.
\end{proof}
\vskip 3mm
In this way, by applying them to the particular log concave function $\exp \left( -\frac{\|\cdot\|_K^2}{2} \right)$, 
many of the known inequalities for convex bodies can be deduced from the corresponding ones for log concave functions. 
Examples include:
\par
\noindent
(i) For $-\infty \leq p \leq \infty$, $p \neq -n$,  the function $f(t)=t^\frac{p}{n+p}$ is concave for $p \geq 0$ and convex for $p \leq 0$. Then,   as a consequence of Theorem \ref{thm00} applied to $\varphi =  \exp \left( -\frac{\|\cdot\|_K^2}{2}\right )$,  and the  Blaschke Santal\'o inequality \cite{Blaschke, MeyerPajor90, Santalo} we obtain the $L_p$-affine isoperimetric inequalities of  \cite{Lutwak1996, WernerYe2008}.  If we apply Proposition  \ref{theo:mono} to  $\varphi =  \exp \left( -\frac{\|\cdot\|_K^2}{2} \right)$, we obtain the  monotonicity behavior of the $L_p$-affine surface area for convex bodies proved in  \cite{Lutwak1996, WernerYe2008}.
\vskip 2mm
\noindent
(ii) If we apply   Proposition \ref{affineinvariant-prop}  and Corollary \ref{affineinvariant-cor} to $\varphi = \exp \left( -\frac{\|\cdot\|_K^2}{2} \right)$,  then we obtain
entropy inequalities first proved in \cite{PaourisWerner2011} for convex bodies. E.g, Corollary \ref{affineinvariant-cor}, together with the Blaschke Santal\'o inequality gives 
the following isoperimetric inequality of \cite{PaourisWerner2011}
\begin{equation*}\label{iso-inequal1}
\Omega_{K^\circ} \ \leq \ \Omega_{\bigg( \frac{B^n_2}{ |B^n_2 |^{\frac{1}{n}}} \bigg)^\circ } \ = \ |B^n_2 |^{2n}.
\end{equation*}
\vskip 3mm
Finally, the duality formula (\ref{polar-identity}) applied to $\varphi =  \exp \left( -\frac{\|\cdot\|_K^2}{2} \right)$ corresponds to a duality formula for convex bodies 
and is a generalization of previously established duality formulas \cite{Hug1996,  WernerYe2008} for convex bodies,  $ as_p(K)=as_{\frac{n^2}{p}}(K^\circ)$. See also \cite{Ludwig2010}.
We skip the proof. 
\vskip 2mm
\begin{prop}\label{}
Let $K$ be a convex body in $\mathbb{R}^n$ and let  $f:(0,  \infty) \rightarrow \mathbb{R}$ be a convex or concave  function. Then 
$$
D_f(P_{K^\circ}, Q_{K^\circ}) = D_{f^*}(P_{K}, Q_{K}).
$$
\end{prop}

\vskip 4mm
\section{Linearization.} \label{Section-lin}

In this section we linearize the  inequalities of Corollary \ref{cor-mono}  around its equality case. 
We treat  only one inequality. The other one is done in the same way.  
We rewrite the 
inequality in terms of a convex function
$\psi: \mathbb R^{n} \to \mathbb R$ such that
 $\varphi=e^{-\psi}$ and get
\begin{align}\label{mainineq55}
\int_{\mathbb R^{n}} e^{-\psi} f \left( e^{2 \psi - \langle \nabla \psi, x \rangle}\det(\text{Hess}(\psi))\right) dx  \geq 
  f \left(\frac{(2\pi)^n}{(\int_{\mathbb R^{n}} e^{-\psi}dx)^2} \right) 
\left(\int_{\mathbb R^{n}} e^{-\psi}dx\right).
\end{align}
Corollary \ref{cor-mono}   requires that $\varphi$ has center of mass at the origin. This is the case if $\psi$ is even. 
We then linearize
around the equality case $\psi(x) =   \|x\|^2/2$ and
obtain the following  functional inequalities. See also \cite{ArtKlarSchuWer},  \cite{HoudreKagan1995}, \cite{HoudrePerez1995}.
The proof, which we include for completeness,  follows \cite{ArtKlarSchuWer}.
Throughout,  
$\|  . \|_{HS}$ denotes the Hilbert Schmidt norm and $\triangle \eta = \operatorname{tr}( \text{Hess}\  \eta)$
is the Laplacian of $\eta$. 
$\gamma_n$ is the normalized Gaussian measure on $\mathbb{R}^n$ and 
$\mbox{Var}_{\gamma_n}(\eta)= \int_{\R^{n}} \eta^2
d\gamma_n - \left(\int_{\R^{n}} \eta  d\gamma_n\right)^2$
is the variance of $\eta$. 
\vskip 3mm
\begin{cor}\label{Cor}
Let $\eta\in C^{2}( \R^{n})\cap
L^{2}(\R^{n}, \gamma_{n})$ be even.  Then
\par
\noindent
(i) 
 \hskip 5mm $ \frac{1}{2} \int_{\mathbb R^{n}} \left(\triangle \eta - \langle \nabla \eta, x \rangle \right )^2 \ d\gamma_n  \ \leq \
 \int_{\mathbb R^{n}} \| \text{Hess} \
\eta \|_{HS}^2  \ d\gamma_n.$
\vskip 2mm
\noindent
\begin{eqnarray*}
&& \hskip -13mm (ii) \hskip 5mm  \int_{\mathbb R^{n}}  \| \nabla \eta \|^2  \ d\gamma_n - \frac{1}{4}  \int_{\mathbb R^{n}} \left(\triangle \eta - \langle \nabla \eta, x \rangle \right )^2 \ d\gamma_n \leq  \mbox{Var}_{\gamma_n}(\eta) \leq \nonumber  \\
&&\hskip 10mm \int_{\mathbb R^{n}}  \| \nabla \eta \|^2  \ d\gamma_n - \frac{1}{2}  \int_{\mathbb R^{n}} \left(\triangle \eta - \langle \nabla \eta, x \rangle \right )^2 \ d\gamma_n + \frac{1}{2}  \int_{\mathbb R^{n}}  \| \text{Hess} 
\eta \|_{HS}^2 d\gamma_n.
\end{eqnarray*}
\end{cor}
\vskip 2mm
\noindent
{\bf Remark.} 
The left hand side of Corollary \ref{Cor} (ii) together with Corollary \ref{Cor}  (i) gives the following  reverse Poincar\'e inequality obtained  in \cite{ArtKlarSchuWer}  
(see also  \cite{HoudreKagan1995}, \cite{HoudrePerez1995})
\begin{equation*}\label{Poincare}
\int_{\R^{n}}  \left[ \left\|\grad \eta\right\|^2_{} - \frac{ \| \text{Hess} \
\eta \|_{HS}^2}{2}
\right] d\gamma_n \le \mbox{Var}_{\gamma_n}(\eta).
\end{equation*}
\vskip 2mm

\begin{proof}
We first prove the corollary for functions with bounded
support.
Thus, let $\eta$ be an even,  twice continuously differentiable
function with bounded support and let $\psi(x)
= \|x\|^2/2 + \eps \eta(x)$. Note that for sufficiently small
$\eps$ the function $\psi$ is convex and that $\psi$ is even as $\eta$ is even.Therefore we can plug $\psi$ into inequality
(\ref{mainineq55}) and develop in
powers of $\eps$. We evaluate first the right hand
expression of (\ref{mainineq55}). 
\begin{eqnarray*}
&f&\left(\frac{(2\pi)^n}{\left(\int_{\mathbb R^{n}} e^{-\|x\|^2/2 - \eps \eta}\ dx\right)^2} \right) 
\left(\int_{\mathbb R^{n}} e^{-\|x\|^2/2 - \eps \eta}\ dx \right)  \\
&=&   (2\pi)^{n/2} \  f \left(1 +2  \eps \int_{\mathbb R^{n}} \eta d\gamma_n  - \eps^2 \int_{\mathbb R^{n}} \eta^2 d\gamma_n +3 \eps^2 
\left( \int_{\mathbb R^{n}} \eta d\gamma_n\right)^2 \right) \\ 
&&\hskip 40mm \bigg(1 -  \eps \int_{\mathbb R^{n}} \eta d\gamma_n 
 +  \frac{\eps^2}{2} \int_{\mathbb R^{n}} \eta^2 d\gamma_n \bigg) + O(\eps^3). \\
\end{eqnarray*}
As $f(1+t)=f(1) +f^\prime(1) t + \frac{f^{\prime \prime}(1)}{2}  t^2 +O(t^3)$, we get that the right hand side of (\ref{mainineq55}) equals
\begin{eqnarray*}
&&  (2\pi)^{n/2} \ \bigg( f (1) + \eps \left[2f^\prime(1)-f(1)\right] \  \int_{\mathbb R^{n}} \eta d\gamma_n + \eps^2 \bigg[ \int_{\mathbb R^{n}} \eta^2 d\gamma_n\left(\frac{f(1)}{2} - f^\prime(1)\right) +\\
&& \left(\int_{\mathbb R^{n}} \eta d\gamma_n\right)^2 \left[ f^\prime(1) +2 f^{\prime \prime}(1)\right]\bigg]\bigg) +
O(\eps^3).
\end{eqnarray*}
\par
We evaluate now the left hand
expression of (\ref{mainineq55}).
Since
$
\text{Hess}\ \psi=I+\eps \ \text{Hess}\  \eta
$,
we obtain for the left hand side
\begin{equation*}
 \int_{{\mathbb R}^n} e^{-\|x\|^2/2 - \eps \eta} f\left(e^{\eps(2 \eta -\langle \nabla \eta, x \rangle)} \det \left(I
+ \eps \ \text{Hess} \  \eta\right)\right) dx.
\end{equation*}
By Taylor's theorem this equals
\begin{eqnarray*}
&& \hskip -10mm \int_{{\mathbb R}^n} e^{-\|x\|^2/2} \left(1 - \eps \eta +
\frac{\eps^2}{2} \eta^2 \right) \cdot
 f \bigg( \left(1 + \eps (2 \eta -\langle \nabla \eta, x \rangle)  +
\frac{\eps^2}{2} (2 \eta -\langle \nabla \eta, x \rangle)^2 \right) \\
&&\hskip 40mm \cdot \det (I
+ \eps  \ \text{Hess}\  \eta)\bigg) dx
+ O(\eps^3).
\end{eqnarray*}
For a matrix $A=(a_{i,j})_{i,j=1,\ldots,n}$, let $D(A) = \sum_{i
=1}^n\sum_{j\neq i}^n [a_{i,i}a_{j,j} - a_{i,j}^2]$. Note that each
$2\times 2$ minor is counted twice. Then
$
\det(I+\eps \ \text{Hess}\  \eta)
=1+\eps \triangle \eta
+\frac{\eps^{2}}{2}D(\text{Hess}\  \eta)+O(\eps^{3})
$.
Therefore
the left hand side of (\ref{mainineq55}) equals
\begin{eqnarray*}
&&  (2\pi)^{n/2} \ \bigg[\int_{\mathbb R^{n}} \left(1 - \eps \eta +
\frac{\eps^2}{2} \eta^2 \right) f\bigg(1 + \eps \left( 2 \eta + \triangle \eta - \langle \nabla \eta, x \rangle \right) \\
&&+ \eps^2 \left( \frac{\left( 2 \eta - \langle \nabla \eta, x \rangle\right)^2}{2} + \triangle \eta \left(2 \eta -  \langle \nabla \eta, x \rangle\right) +
\frac{D(\text{Hess}\  \eta)}{2}  \right) \bigg)
   d\gamma_n 
\bigg]+O(\eps^{3}) = \\
&&(2\pi)^{n/2} \ \bigg( f(1) + \eps \left[  -f(1)\int_{\mathbb R^{n}} \eta d\gamma_n + f^\prime(1) \left(\int_{\mathbb R^{n}}(2 \eta-  \langle \nabla \eta, x \rangle + \triangle \eta)d\gamma_n \right)  \right]  \\
&&+ \eps^2\bigg[f(1) \int_{\mathbb R^{n}} \frac{\eta^2}{2} d\gamma_n - f^\prime(1) \bigg( \int_{\mathbb R^{n}} \eta \left( 2 \eta - \langle \nabla \eta, x \rangle + \triangle \eta \right)  d\gamma_n\bigg) + \\
&& f^\prime(1) \left(\int_{\mathbb R^{n}}\left(\frac{\left( 2 \eta - \langle \nabla \eta, x \rangle\right)^2}{2} + \triangle \eta \left(2 \eta -  \langle \nabla \eta, x \rangle\right) +
\frac{D(\text{Hess}\  \eta)}{2} \right) d \gamma_n\right) \\
&&+ \frac{f^{\prime \prime} (1)}{2} \left( \int_{\mathbb R^{n}} \left( 2 \eta-  \langle \nabla \eta, x \rangle + \triangle \eta\right)^2d \gamma_n\right) \ \bigg] \bigg)+O(\eps^{3}).
\end{eqnarray*}
Now observe that   $\int_{\mathbb R^{n}} (\triangle \eta -  \langle \nabla \eta, x \rangle )d\gamma_n=0$.  Also,  as  the coefficients of  order zero  and  of order $\eps$ are the same on the left hand side  and the right hand side,  we discard them. We divide both sides by $\eps^2$ and take the limit for $\eps \rightarrow 0$.  Thus,  the inequality (\ref{mainineq55}) is equivalent to 
\begin{eqnarray*}
&&\hskip -5mm \bigg( f^\prime(1) + 2 f^{\prime \prime}(1) \bigg) \left[ \left(\int_{\mathbb R^{n}} \eta d\gamma_n\right)^2  d\gamma_n - \int_{\mathbb R^{n}} \eta^2  \ d\gamma_n \right]   \leq \\
&& \hskip -5mm \bigg( f^\prime(1) + 2 f^{\prime \prime}(1) \bigg) \left[  \int_{\mathbb R^{n}} \eta  \triangle \eta \ d\gamma_n - \int_{\mathbb R^{n}}   \eta \langle \nabla \eta, x \rangle \ d\gamma_n \right]   + 
 \frac{f^\prime(1)}{2} \int_{\mathbb R^{n}}  D(\text{Hess}\ \eta) \ d\gamma_n    \\
&&\hskip -5mm +   \frac{f^{\prime \prime}(1)}{2} \int_{\mathbb R^{n}} (\triangle \eta)^2 d\gamma_n  
  + \bigg( f^\prime(1) +  f^{\prime \prime}(1) \bigg) \left[ \int_{\mathbb R^{n}}  \frac{\langle \nabla \eta, x \rangle^2}{2}  d\gamma_n   - \int_{\mathbb R^{n}}  \triangle \eta \langle \nabla \eta, x \rangle  d\gamma_n  \right].   
\end{eqnarray*}
\noindent
Integration by parts yields
$
\int_{\mathbb R^{n}}   \eta \langle \nabla \eta, x \rangle  d\gamma_n \ = \ \frac{1}{2}   \int_{\mathbb R^{n}}   \eta^2 (x) \ ( \|x\|^2  - n ) \ d\gamma_n 
$
\noindent
and
$$
 \int_{\mathbb R^{n}} \eta  \triangle \eta \ d\gamma_n =  -   \int_{\mathbb R^{n}}  \| \nabla \eta \|^2  \ d\gamma_n  + \frac{1}{2}   \int_{\mathbb R^{n}}   \eta^2 (x) \ ( \|x\|^2  - n ) \ d\gamma_n  .
$$
We put $a=f^\prime(1)$ and  $ b= f^{\prime \prime}(1) $.
Note that $a \leq 0$ and $ b \geq 0$,  as $f$ is   convex and  decreasing.
Thus,  the inequality becomes
\begin{eqnarray*}
&& \hskip -10mm \big( a + 2b \big)  \left( \mbox{Var}_{\gamma_n}(\eta) -  \int_{\mathbb R^{n}}  \| \nabla \eta \|^2  \ d\gamma_n  \right) \nonumber \\   
&& \geq  \frac{a}{2} \int_{\mathbb R^{n}} \| \text{Hess} \
\eta \|_{HS}^2  \ d\gamma_n     -   \frac{a+b}{2} \int_{\mathbb R^{n}} \left(\triangle \eta - \langle \nabla \eta, x \rangle \right )^2 \ d\gamma_n . 
\end{eqnarray*}
\par
Hence we have shown that the inequality holds for all twice continuously
differentiable functions $\eta$ with bounded support.
Now we extend it to
all twice continuously differentiable functions $\eta\in L^{2}(\mathbb R^{n},\gamma_{n})$
satisfying  the necessary integrability conditions, 
by a standard approximation argument, as follows.
\par
Let $\chi_{k}$ be a twice continuously differentiable function
bounded between zero and one such that $\chi_{n}(x)=1$ for all
$\|x\|\leq k$ and $\chi_{n}(x)=0$ for all $\|x\|>k+1$. Then, for all
$k\in\mathbb N$
\begin{eqnarray*}
&& \hskip -10mm \big( a + 2b \big)  \left( \mbox{Var}_{\gamma_n}(\eta\cdot \chi_{k}) -  \int_{\mathbb R^{n}}   \| \nabla ( \eta\cdot \chi_{k} ) \|^2   \ d\gamma_n  \right) \nonumber \\   
&& \geq  \frac{a}{2} \int_{\mathbb R^{n}} \| \text{Hess} \
(\eta\cdot \chi_{k}) \|_{HS}^2  \ d\gamma_n     -   \frac{a+b}{2} \int_{\mathbb R^{n}} \bigg( \triangle (\eta\cdot \chi_{k}) - \langle \nabla ( \eta\cdot \chi_{k}), x \rangle  \bigg)^2 \ d\gamma_n ,
\end{eqnarray*}
which is equivalent to
\begin{eqnarray*}
&&       2b \left[ \left(\int_{\R^{n}} ( \eta\cdot\chi_{k} ) d\gamma_n\right)^2 + \int_{\mathbb R^{n}}   \| \nabla ( \eta\cdot \chi_{k} ) \|^2   \ d\gamma_n    \right]         \\
& - & \ a \left[  \int_{\R^{n}} ( \eta\cdot\chi_{k} )^2 
d\gamma_n \ + \ \int_{\mathbb R^{n}} \frac{1}{2}\bigg( \triangle (\eta\cdot \chi_{k}) - \langle \nabla ( \eta\cdot \chi_{k}), x \rangle  \bigg)^2\ d\gamma_n     \right] \\
&  \leq & ( -a ) \left[ \left(\int_{\R^{n}} ( \eta\cdot\chi_{k} ) d\gamma_n\right)^2 +  \int_{\mathbb R^{n}}   \| \nabla ( \eta\cdot \chi_{k} ) \|^2   \ d\gamma_n  +  \int_{\mathbb R^{n}} \frac{\| \text{Hess}  \
\eta \|_{HS}^2 }{2} \ d\gamma_n  \right]   \nonumber \\ 
& + & \ b \left[  \int_{\R^{n}} 2( \eta\cdot\chi_{k} )^2
d\gamma_n \ + \ \int_{\mathbb R^{n}} \frac{1}{2}\bigg( \triangle (\eta\cdot \chi_{k}) - \langle \nabla ( \eta\cdot \chi_{k}), x \rangle  \bigg)^2\ d\gamma_n    \right].
\end{eqnarray*}
Now we pass to the limit $k \rightarrow \infty$ on both sides and obtain
\begin{eqnarray*}
&&  2b \left[ \liminf_{k\to\infty} \left( \int_{\R^{n}} ( \eta\cdot\chi_{k} ) d\gamma_n\right)^2 + \liminf_{k\to\infty}   \int_{\mathbb R^{n}}   \| \nabla ( \eta\cdot \chi_{k} ) \|^2   \ d\gamma_n    \right]         \\
& - & a \left[ \liminf_{k\to\infty}  \int_{\R^{n}} ( \eta\cdot\chi_{k} )^2 
d\gamma_n \ +  \liminf_{k\to\infty}  \int_{\mathbb R^{n}} \frac{1}{2}\bigg( \triangle (\eta\cdot \chi_{k}) - \langle \nabla ( \eta\cdot \chi_{k}), x \rangle  \bigg)^2\ d\gamma_n     \right] \\
& \leq &  (-a)\left[\limsup_{k\to\infty} \left(\int_{\mathbb R^{n}}  ( \eta\cdot\chi_{k} ) d\gamma_n\right)^2 + \limsup_{k\to\infty} \int_{\mathbb R^{n}}  \| \nabla ( \eta\cdot \chi_{k} ) \|^2    d\gamma_n   \right]    +\nonumber \\
&b& \left[ \limsup_{k\to\infty} \int_{\R^{n}} 2( \eta\cdot\chi_{k} )^2
d\gamma_n  +  \limsup_{k\to\infty} \int_{\mathbb R^{n}} \frac{1}{2}\bigg( \triangle (\eta\cdot \chi_{k}) - \langle \nabla ( \eta\cdot \chi_{k}), x \rangle  \bigg)^2\ d\gamma_n    \right] \\
&-&  a \left[  \limsup_{k\to\infty} \int_{\mathbb R^{n}}  \frac{\| \text{Hess}  \
\eta \|_{HS}^2 }{2} \ d\gamma_n   \right] .
\end{eqnarray*}
Fatou's lemma and the dominated convergence theorem yield
\begin{eqnarray*}
&& \hskip -15mm \frac{a}{2} \int_{\mathbb R^{n}} \| \text{Hess} \
\eta \|_{HS}^2  \ d\gamma_n     -   \frac{a+b}{2} \int_{\mathbb R^{n}} \left(\triangle \eta - \langle \nabla \eta, x \rangle \right )^2 \ d\gamma_n \\
&& \leq \  \big( a + 2b \big)  \left( \mbox{Var}_{\gamma_n}(\eta) -  \int_{\mathbb R^{n}}  \| \nabla \eta \|^2  \ d\gamma_n  \right) . 
\end{eqnarray*}
Finally, we consider the cases $a+2b=0$, $a+2b>0$ and $a+2b <0$ and optimize in each case. 
\end{proof}

\vskip 2mm 
\noindent 
Umut Caglar\\
{\small Department of Mathematics \\
{\small Case Western Reserve University \\
{\small Cleveland, Ohio 44106, U. S. A. \\
{\small \tt umut.caglar@case.edu}\\ \\
\noindent
\and 
Elisabeth Werner\\
{\small Department of Mathematics \ \ \ \ \ \ \ \ \ \ \ \ \ \ \ \ \ \ \ Universit\'{e} de Lille 1}\\
{\small Case Western Reserve University \ \ \ \ \ \ \ \ \ \ \ \ \ UFR de Math\'{e}matique }\\
{\small Cleveland, Ohio 44106, U. S. A. \ \ \ \ \ \ \ \ \ \ \ \ \ \ \ 59655 Villeneuve d'Ascq, France}\\
{\small \tt elisabeth.werner@case.edu}\\ \\

\end{document}